\documentclass[12pt]{article}
\usepackage[T1]{fontenc}
\usepackage{lmodern}
\usepackage{amsmath,amssymb,amsthm,mathtools}
\usepackage{microtype}
\usepackage[top=27mm, bottom=27mm, left=22mm, right=22mm]{geometry}
\usepackage{xcolor}
\usepackage{tikz}
\usepackage{needspace}
\usepackage{placeins}
\usepackage{cite}
\usepackage{array,booktabs,longtable}
\usepackage[colorlinks=true,linkcolor=blue!50!black,citecolor=blue!50!black,urlcolor=blue!50!black]{hyperref}

\newtheorem{theorem}{Theorem}[section]

\newtheorem{corollary}[theorem]{Corollary}
\newtheorem{lemma}[theorem]{Lemma}
\newtheorem{problem}[theorem]{Problem}

\newtheorem{conjecture}{Conjecture}
\newtheorem{claim}[theorem]{Claim}
\theoremstyle{definition}
\newtheorem{definition}[theorem]{Definition}

\newtheorem{remark}[theorem]{Remark}
\newtheorem*{remark*}{Remark}

\newcommand{\F}{\mathbb F}
\newcommand{\Z}{\mathbb Z}
\newcommand{\E}{\mathbb E}
\newcommand{\Prb}{\mathbb P}
\newcommand{\cS}{\mathcal S}
\newcommand{\cA}{\mathcal A}

\newcommand{\cF}{\mathcal F}

\newcommand{\cP}{\mathcal P}

\newcommand{\bfone}{\boldsymbol 1}
\newcommand{\bfzero}{\boldsymbol 0}

\DeclareMathOperator{\GL}{GL}

\title{The probabilistic bound for perfect hashing is sharp only for binary alphabets}

\author{Shitao Li\thanks{School of Internet, Anhui University, Hefei, Anhui 230039, China. E-mail: \textsf{lishitao0216@163.com}.},~ 
Minjia Shi\thanks{Key Laboratory of Intelligent Computing and Signal Processing, Ministry of Education, School of Mathematical Sciences, Anhui University, Hefei, Anhui 230601, China. E-mail: \textsf{smjwcl.good@163.com}}
\thanks{Corresponding author: Minjia Shi.}}
\date{}

\begin{document}
\maketitle

\begin{abstract}
We investigate the optimal asymptotic rate $R_q$ of perfect $q$-hash codes, a fundamental extremal problem in combinatorics, information theory, and computer science. Let $B_q=-\log_2(1-q!/q^q)/(q-1)$ denote the classical probabilistic lower bound. We prove that $R_q=B_q$ if and only if $q=2$, establishing an explicit strict improvement for every $q\ge3$ and settling the first assertion of a conjecture of Xing and Yuan (Combinatorica 43 (2023), 347--366). More significantly, we show that, for every $\eta>0$, $R_q\ge B_q+\exp(-(2+\eta)q)$ for all sufficiently large integers $q$. This improves the previously known increments of scale $\exp(-3q+O(\ln q))$ by an exponential factor. We also establish a structural limitation of the concatenation method. For pairwise uniform inner codes of subexponential length, the resulting rate bound exceeds $B_q$ by at most $\exp(-3q+o(q))$. Thus, exponentially long inner codes are necessary to achieve the new improvement scale within this framework. Our proofs combine Fourier analysis, high-dimensional simplex codes, additive permutation counts, and entropy bounds for hypergraph matchings. We further strengthen lower bounds from short inner codes and improve several known numerical estimates.
\end{abstract}


\section{Introduction}

Let $q\ge2$ be an integer and write $[q]=\{1,\ldots,q\}$. A code $C\subseteq[q]^N$ is a \emph{perfect $q$-hash code} if every $q$ distinct codewords $c_1,\ldots,c_q\in C$ have a coordinate $i\in[N]$ such that $\{c_{1,i},\ldots,c_{q,i}\}=[q]$. Such a coordinate \emph{separates} the codewords. Perfect $3$-hash codes are also called \emph{trifferent codes}. Equivalently, the coordinate projections of a perfect $q$-hash code $C$ form a family of functions from $C$ to $[q]$ such that every $q$-element subset of $C$ is mapped injectively by at least one of them. This formulation as a perfect hash family underlies applications in combinatorics and computer science \cite{FK1984}. For example, such families are used to derandomize the color-coding method of Alon, Yuster and Zwick \cite{AYZ1995}. Explicit perfect hash families with more general alphabet and separation parameters can be constructed from algebraic curves \cite{WX2001}. Perfect $q$-hash codes also describe zero-error list decoding with list size at most $q-1$ on the $q/(q-1)$ channel \cite[Section 1]{XY2023}. In this channel, each input symbol can be received as any of the other $q-1$ symbols, and the separation condition guarantees that no received word is compatible with $q$ distinct codewords.

The basic extremal problem is to determine the largest size $M(N,q)$ of a perfect $q$-hash code of length $N$. For fixed $q$, we study the asymptotic rate
\[ R_q=\limsup_{N\to\infty}\frac{\log_2 M(N,q)}{N},\]
which measures the exponential growth of $M(N,q)$. The study of $R_q$ goes back to the 1980s, but its exact value is not known for any $q\ge3$.

\subsection{Bounds and the Xing--Yuan conjecture}
We first recall the upper bounds. Fredman and Koml\'os \cite{FK1984} proved the general upper bound $R_q\le q!/q^{q-1}$. K\"orner \cite{Korner1986} developed an information-theoretic approach to these bounds, and Nilli \cite{Nilli1994} found a short probabilistic proof of upper bounds for the more general separating-family problem. The Fredman--Koml\'os bound was later improved for $q=4$ by Dalai, Guruswami and Radhakrishnan \cite{DGR2020}, and for every $q\ge5$ by Guruswami and Riazanov \cite{GR2022}.
Costa and Dalai \cite{CD2021} resolved the polynomial optimization conjecture arising in the argument of \cite{GR2022} and obtained further improvements for $q=5,6$. Della Fiore, Costa and Dalai \cite{DCD2022} extended these methods to $(b,k)$-hashing and strengthened several bounds for small parameters.
For trifferent codes, an elementary argument shows that $M(N,3)\le2(3/2)^N$ \cite[Section 1.1]{BDGP2024}, and Bhandari and Khetan \cite{BK2025} proved the stronger estimate $M(N,3)=O\bigl(N^{-2/5}(3/2)^N\bigr)$. Since the factor $N^{-2/5}$ does not affect the exponential rate, the best known upper bound on $R_3$ remains $\log_2(3/2)$.

Progress on lower bounds has been much more limited. Let $p_q=q!/q^q$ be the probability that $q$ independent uniform symbols are all distinct, and let ${B_q=-\log_2(1-p_q)/(q-1)}$.
By choosing codewords independently and deleting one codeword from each nonseparated $q$-subset,
Fredman and Koml\'os \cite{FK1984} and K\"orner \cite{Korner1986} obtained the classical probabilistic lower bound
\begin{equation}\label{eq:probabilistic-bound}
	R_q\ge B_q\qquad(q\ge2).
\end{equation}
Before the present work, this bound had been improved only for finitely many small alphabet sizes and for all sufficiently large $q$.
K\"orner and Marton \cite{KM1988} obtained the first strict improvement in 1988. For $q=3$, they proved
\begin{equation}\label{eq:KMbound}
	R_3\ge\frac14\log_2\frac95 >\frac12\log_2\frac97=B_3
\end{equation}
by concatenating a short ternary inner code with a randomly chosen outer code. 

Later work on trifferent codes includes linear and explicit constructions. Bishnoi, D'haeseleer, Gijswijt and Potukuchi \cite{BDGP2024} proved that ternary linear trifferent codes are equivalent to minimal codes, and hence to strong blocking sets. They also constructed linear trifferent codes of size at least $(9/5)^{N/4}/3$ for all sufficiently large
$N$ \cite[Theorem 1.5]{BDGP2024}, which matches the exponential lower
bound in \eqref{eq:KMbound} within the class of linear codes. Explicit
constructions were developed further by Alon, Bishnoi, Das and Neri
\cite{ABDN2024}, using expander graphs, and by Bishnoi and Tomon
\cite{BT2026}, using blocking sets of higher codimension. These results concern
linearity and explicitness, rather than a strict improvement on
\eqref{eq:probabilistic-bound} for every $q$. Recently, Han, Sun and Ge \cite[Theorem 1.1]{HSG2026} proved that $M(N,3)\ge c\sqrt N(9/5)^{N/4}$ for some absolute constant $c>0$ and all sufficiently large $N$. This improves the lower bound on $M(N,3)$ by an unbounded factor, but leaves the lower bound on $R_3$ unchanged. Thus, \eqref{eq:KMbound} is still the best known lower bound on $R_3$.

For more than three decades, $q=3$ remained the only alphabet size for which a strict improvement on \eqref{eq:probabilistic-bound} was known \cite{XY2023}. The next advance is due to Xing and Yuan, who removed the requirement that the inner code itself be a perfect hash code. Instead, they count the separated $q$-element subsets of the inner code and choose an outer code adapted to this family. This additional freedom permits maximum distance separable (MDS) inner codes over finite abelian groups and relates the rate problem to the enumeration of simultaneous permutations. Using the asymptotic counts of additive triples of bijections due to Eberhard, Manners and Mrazovi\'c \cite{EMM2019} and Eberhard \cite{Eberhard2017}, Xing and Yuan proved
\begin{align}
R_q&\ge-\frac1{4(q-1)}\log_2\left((1-p_q)^4-\left(\frac{3q}{\sqrt e}+o(q)\right)p_q^3\right),
&q\not\equiv2\pmod4,\label{eq:xy-regular}\\
R_q&\ge-\frac1{3(q-1)}\log_2\left((1-p_q)^3-\left(\frac{q}{2\sqrt e}+o(q)\right)p_q^3\right),
&q\equiv2\pmod4.\label{eq:xy-exceptional}
\end{align}
These bounds exceed $B_q$ for all sufficiently large $q$ \cite[Theorem 1.1]{XY2023}. Xing and Yuan \cite{XY2023} also obtained strict improvements for $4\le q\le15$ and for every odd $q$ with $17\le q\le25$, and reported numerical evidence for several further odd values of $q$. To the best of our knowledge, \eqref{eq:xy-regular} and \eqref{eq:xy-exceptional} remain the best known general asymptotic lower bounds on $R_q$ as $q\to\infty$.

These results left open whether the probabilistic bound can be improved for every alphabet size, and not only for the small values verified in \cite{XY2023} and for all sufficiently large values. Xing and Yuan formulated a conjecture consisting of two assertions. 

\begin{conjecture}[Xing--Yuan {\cite[Conjecture 1.2]{XY2023}}]
For every integer $q\ge3$, there exist families of perfect $q$-hash codes whose asymptotic rate is strictly larger than $B_q$. Moreover, this strict improvement can be obtained from the concatenation constructions specified in Theorems~4.7, 4.8 and~5.2 of \cite{XY2023}, in their respective ranges of alphabet sizes.
\end{conjecture}

In this paper, we prove the first assertion of the Xing--Yuan conjecture with an explicit lower bound for every $q\ge3$. Once an improvement is known to exist for every alphabet size, the next question is how large it can be. We improve \eqref{eq:xy-regular} and \eqref{eq:xy-exceptional} on a different exponential scale, uniformly over all sufficiently large integers $q$, and we show that inner codes of the type used by Xing and Yuan cannot reach this scale unless their length grows exponentially with $q$. 

\subsection{Main results}

Our first theorem answers the qualitative question. It characterizes the alphabet sizes for which the probabilistic lower bound \eqref{eq:probabilistic-bound} can be improved.

\begin{theorem}\label{thm:allq}
	One has $R_2=B_2=1$. For every integer $q\ge3$,
	\begin{equation}\label{eq:allq}
		R_q\ge -\frac{1}{(2q+1)(q-1)}
		\log_2\left((1-p_q)^{2q+1}-p_q^{2q+1}\right)>B_q.
	\end{equation}
	In particular,
	\begin{equation}\label{eq:sharpness}
		R_q=B_q\quad\text{if and only if}\quad q=2.
	\end{equation}
\end{theorem}

Theorem~\ref{thm:allq} proves the first assertion of the Xing--Yuan conjecture for every $q\ge3$, with no unspecified threshold and no finite range of $q$ left to verification. The lower bound in \eqref{eq:allq} is obtained by concatenating an explicit inner code of length $2q+1$ with a randomly chosen outer code. Theorem~\ref{thm:length} replaces this length by a computable odd integer $\ell_*(q)\le(1+o(1))q/\ln q$. The equivalence \eqref{eq:sharpness} characterizes equality in the probabilistic lower bound and does not determine $R_q$.

We next ask how large the improvement can be. By Stirling's formula, we have
\begin{equation}\label{eq:scale}
	p_q\sim\sqrt{2\pi q}\,e^{-q}\quad\text{and}\quad B_q\sim\frac{p_q}{(q-1)\ln2}.
\end{equation}
Expanding the logarithms in \eqref{eq:xy-regular} and \eqref{eq:xy-exceptional} as in Section~\ref{sec:concat}, we see that these bounds exceed $B_q$ by $(3/(4\sqrt e\ln2)+o(1))p_q^3$ and $(1/(6\sqrt e\ln2)+o(1))p_q^3$, respectively. Both increments have the exponential scale $\exp(-3q+O(\ln q))$. Our main asymptotic result provides a much larger increment over the probabilistic lower bound, uniformly over all sufficiently large alphabet sizes.

\begin{theorem}\label{thm:exponential}
	For every real $\eta>0$, there exists an integer $q_0(\eta)$ such that
	\begin{equation}\label{eq:exponential}
		R_q\ge B_q+\exp\bigl(-(2+\eta)q\bigr) \quad\text{for every integer }q\ge q_0(\eta).
	\end{equation}
	Consequently, we have
	\begin{equation}\label{eq:liminf}
		\liminf_{q\to\infty}\frac{1}{q}\ln(R_q-B_q)\ge-2.
	\end{equation}
\end{theorem}

In contrast with the increments of Xing and Yuan, Theorem~\ref{thm:exponential} implies
\begin{equation}\label{eq:superpolynomial}
	\frac{R_q-B_q}{q^A p_q^3}\longrightarrow\infty
	\quad\text{for every fixed }A>0,
\end{equation}
where $q$ ranges over all integers. These improvements concern only the excess over $B_q$, and the leading order of the lower bound is still $B_q$. Theorem~\ref{thm:exponential} does not determine the exponential order of $R_q-B_q$, and unlike Theorem~\ref{thm:allq}, it does not specify the threshold $q_0(\eta)$.

We prove Theorem~\ref{thm:exponential} by establishing a more general quantitative bound for a family of high-dimensional simplex inner codes. This bound allows the inner length to grow exponentially with the alphabet size and also yields polynomial improvements when the dimension is fixed. Let $s\ge3$ be a prime-power divisor of $q$, let $r=q/s$, and let $\cP_{s,d}$ contain one nonzero vector from each one-dimensional subspace of $\F_s^d$. With $H=\Z/r\Z$, we use the alphabet $\F_s\times H$ of size $q$ and define
\begin{equation}\label{eq:simplex-code}
	C(q,s,d)=\bigl\{\bigl((u\cdot v,b_v)\bigr)_{v\in\cP_{s,d}}:u\in\F_s^d,\ (b_v)_{v\in\cP_{s,d}}\in H^\ell\bigr\},\qquad \ell=\frac{s^d-1}{s-1}.
\end{equation}
The first components form the $s$-ary simplex code of dimension $d$, and the second components are unrestricted, so that $|C(q,s,d)|=s^dr^\ell$. If $r=1$, then $H$ is trivial.

\begin{theorem}\label{thm:general}
	Fix $0<\varepsilon<1$. Let $s\ge3$ be a prime-power divisor of $q$, let $r=q/s$, and let $d\ge2$ satisfy
	\begin{equation}\label{eq:lengthcondition}
		\ell=\frac{s^d-1}{s-1}\le\exp((1-\varepsilon)q).
	\end{equation}
	Set $a_{q,s}=se^{-1/(2r)}-1$. Then, as $q\to\infty$, concatenation with the inner code \eqref{eq:simplex-code} yields
	\begin{equation}\label{eq:generalrate}
		R_q\ge B_q+(1-o(1)) \frac{(\ell-1)(s-1)a_{q,s}}{6(q-1)\ln2}\,p_q^3,
	\end{equation}
	where the $o(1)$ term depends only on $\varepsilon$ and $q$. In particular, it is uniform over all such $s$ and $d$.
\end{theorem}

Since the increment in \eqref{eq:generalrate} grows linearly with $\ell$, Theorem~\ref{thm:exponential} follows by taking $d$ as large as \eqref{eq:lengthcondition} permits. For each fixed $q$, all parameters of the inner code are fixed before the outer length tends to infinity, so allowing $d$ to grow with $q$ in Theorem~\ref{thm:general} does not interchange these two limits. For fixed dimensions, we obtain the following polynomial improvements over the increments of Xing and Yuan.

\begin{corollary}\label{cor:fixedd}
	The following asymptotic bounds hold.
	\begin{enumerate}
		\item [(1)] For each fixed integer $d\ge2$, as $q\to\infty$
		through prime powers,
		\[ R_q\ge B_q+ \left(\frac{q^d}{6\sqrt e\ln2}+o(q^d)\right)p_q^3.\]
		\item [(2)] For each fixed positive integer $r$, as
		$s\to\infty$ through prime powers and $q=rs$,
		\[ R_q\ge B_q+ \left(\frac{e^{-1/(2r)}}{6r^3\ln2}+o(1)\right)q^2p_q^3.\]
	\end{enumerate}
\end{corollary}

For $d=2$ and $s=q$, the code \eqref{eq:simplex-code} is equivalent to the extended Reed--Solomon code of dimension two and length $q+1$ (Remark~\ref{rem:rs}). For this inner code, the case $d=2$ of Corollary~\ref{cor:fixedd}(1) provides the increment $(q^2/(6\sqrt e\ln2)+o(q^2))p_q^3$.

We now turn to inner codes of length at most $q+1$, which include the inner codes of lengths four and three behind \eqref{eq:xy-regular} and \eqref{eq:xy-exceptional}, respectively. For such lengths, a second construction, based on MacNeish's product construction \cite{MacNeish1922}, is sometimes better than the case $d=2$ of Theorem~\ref{thm:general}. Its length is governed by the following arithmetic function. For an integer $n\ge2$, let $\nu(n)$ be the smallest prime power $p^a$ such that $p^a\mid n$ and $p^{a+1}\nmid n$. For $q\ge3$, let
\[(n_q,\gamma_q)= \begin{cases}
	\bigl(\nu(q),\,e^{-1/2}\bigr), & q\not\equiv2\pmod4,\\[2pt]
	\bigl(\nu(q/2),\,1/(2e^{1/4})\bigr), & q\equiv2\pmod4.
\end{cases}\]
We have $n_q\ge3$ for every $q\ge3$, since, if ${q\not\equiv2\pmod4}$, every prime power exactly dividing $q$ is odd or divisible by $4$, and if ${q\equiv2\pmod4}$, then $q/2$ is odd and at least $3$.

\begin{theorem}\label{thm:asym}
	As $q\to\infty$,
	\[
	R_q\ge B_q+(1+o(1))\,\frac{\gamma_q\,n_q(n_q-1)}{6\ln2}\,p_q^3.
	\]
\end{theorem}

The inner code behind Theorem~\ref{thm:asym} has length $n_q+1$. If ${q\not\equiv2\pmod4}$, it is a two-dimensional linear code over the ring $\prod_i\F_{q_i}$, where $q=\prod_iq_i$ is the factorization of $q$ into prime powers, and any two of its coordinates determine a codeword. If ${q\equiv2\pmod4}$, the inner code is obtained from such a code over a ring of order $q/2$ by adjoining free binary components.

\begin{remark*}
Let $R$ be a commutative ring of order $q$, and let $M$ be a matrix over $R$ with two rows such that any two coordinates of $(x,y)M$ determine $(x,y)\in R^2$. If ${q\not\equiv2\pmod4}$, then the code $\{(x,y)M:x,y\in R\}$ has length at most $n_q+1$ by Remark~\ref{rem:ringlength}, so no such code is longer than the inner code behind Theorem~\ref{thm:asym}. For every $q\ge3$, the constant $\gamma_q$ comes from the projections of the inner code behind Theorem~\ref{thm:asym} onto three coordinates. By Remark~\ref{rem:cosetoptimal}, $\gamma_q$ is asymptotically the largest constant obtainable from Lemma~\ref{lem:cosettriples} for three-coordinate projections of the form $\{(x,y,z):x+y+z\in U\}$ with $U$ a subgroup.
\end{remark*}

Since $n_q\ge3$, Theorem~\ref{thm:asym} contains the following bounds.

\begin{corollary}\label{cor:special}
	The following lower bounds hold as $q\to\infty$ through the indicated sets of integers.
	\begin{enumerate}
		\item [(1)] If $q\not\equiv2\pmod4$, then
		\[ R_q\ge B_q+\left(\frac{1}{\sqrt e\ln2}+o(1)\right)p_q^3.\]
		\item [(2)] If $q\equiv2\pmod4$, then
		\[ R_q\ge B_q+\left(\frac{1}{2e^{1/4}\ln2}+o(1)\right)p_q^3.\]
	\end{enumerate}
\end{corollary}

\begin{remark*}
Corollary~\ref{cor:special} increases the leading coefficients of the increments in \eqref{eq:xy-regular} and \eqref{eq:xy-exceptional} by the factors $4/3$ and $3e^{1/4}$, respectively. If $n_q\ge4$, which holds unless $3$ exactly divides $q$, then $n_q(n_q-1)/6\ge2$, and these factors at least double. If $q$ is a prime power, then Theorem~\ref{thm:asym} and Corollary~\ref{cor:fixedd}(1) with $d=2$ have the same leading term $q^2p_q^3/(6\sqrt e\ln2)$. Similarly, if ${q\equiv2\pmod4}$ and $q/2$ is a prime power, then Theorem~\ref{thm:asym} and Corollary~\ref{cor:fixedd}(2) with $r=2$ have the same leading term $q^2p_q^3/(48e^{1/4}\ln2)$. In general, neither of the two constructions of length at most $q+1$ dominates the other. If $q=p_1p_2$ for primes $p_1<p_2<2p_1$, then the increments in Theorem~\ref{thm:asym} and in Theorem~\ref{thm:general} with $d=2$ are of orders $qp_q^3$ and $q^{1/2}p_q^3$, respectively. If $q=3\cdot2^a$ with $a\ge2$, these orders are $p_q^3$ and $q^2p_q^3$, respectively.
\end{remark*}

Table~\ref{tab:increments} summarizes the increments over $B_q$ as $q\to\infty$, where $\ell$ denotes the inner length. The first five rows list the leading terms of these increments, up to a factor $1+o(1)$, whereas the last row contains an explicit lower bound valid for all sufficiently large $q$. In Corollary~\ref{cor:fixedd}(1), $d\ge2$ is fixed, and in Corollary~\ref{cor:fixedd}(2), $r$ is a fixed positive integer and $s$ ranges over prime powers. In the last row, $0<\eta<1$ is fixed.

\begin{table}[t]
\centering
\caption{Comparison of asymptotic rate improvements.}
\label{tab:increments}
\small
\renewcommand{\arraystretch}{1.5}
\setlength{\tabcolsep}{20pt}
\begin{tabular}{@{}llll@{}}
\toprule
Result & $q$ & $\ell$ & Increment\\
\midrule
\cite{XY2023}, \eqref{eq:xy-regular}& $q\not\equiv2\pmod4$& $4$& $\displaystyle\frac{3p_q^3}{4\sqrt e\ln2}$\\

\addlinespace[3pt]

\cite{XY2023}, \eqref{eq:xy-exceptional}& $q\equiv2\pmod4$& $3$& $\displaystyle\frac{p_q^3}{6\sqrt e\ln2}$\\

\addlinespace[3pt]

Theorem~\ref{thm:asym}& all $q$& $n_q+1$& $\displaystyle\frac{\gamma_qn_q(n_q-1)p_q^3}{6\ln2}$\\

\addlinespace[3pt]

Corollary~\ref{cor:fixedd}(1)& prime powers& $\displaystyle\frac{q^d-1}{q-1}$& $\displaystyle\frac{q^dp_q^3}{6\sqrt e\ln2}$\\

\addlinespace[3pt]

Corollary~\ref{cor:fixedd}(2)& $q=rs$& $q/r+1$& $\displaystyle\frac{e^{-1/(2r)}q^2p_q^3}
{6r^3\ln2}$\\

\addlinespace[3pt]

Theorem~\ref{thm:exponential}& all $q$& $\le e^{(1-\eta/2)q}$& $\ge e^{-(2+\eta)q}$\\
\bottomrule
\end{tabular}
\end{table}

The increments in Table~\ref{tab:increments} grow with the inner length. The next theorem shows that this growth is necessary for a large class of inner codes. We call a code $C\subseteq\cA^\ell$ over an alphabet $\cA$ of size $q$ \emph{pairwise uniform} if, for any two distinct coordinates, all pairs of symbols have the same number of preimages under the corresponding projection. All inner codes in this paper are pairwise uniform, and so are the inner codes behind \eqref{eq:xy-regular} and \eqref{eq:xy-exceptional}. Let $F(C)$ be the probability that $q$ codewords chosen independently and uniformly from $C$ are not separated by any coordinate. By the concatenation inequality of Xing and Yuan (Lemma~\ref{lem:concatenation}), we have $R_q\ge-\log_2F(C)/(\ell(q-1))$ whenever $0<F(C)<1$.

\begin{theorem}\label{thm:limit}
Let $q\ge3$, and let $C\subseteq\cA^\ell$ be a pairwise uniform code of length $\ell\ge3$ over an alphabet $\cA$ of size $q$. Let $\theta_q=p_q\exp\bigl(-2q+2(3+\ln2)\sqrt q\bigr)$.
If $\ell p_q\le1/2$ and $\ell^3\theta_q\le1$, then
\[ -\frac{\log_2F(C)}{\ell(q-1)}\le B_q+\frac{\ell^2\theta_q}{q-1}.\]
\end{theorem}

By \eqref{eq:scale}, $\theta_q=p_q^3\exp(O(\sqrt q))$. Hence, uniformly over pairwise uniform inner codes of length $\ell=e^{o(q)}$, the lower bound from Lemma~\ref{lem:concatenation} is at most $B_q+p_q^{3+o(1)}$ as $q\to\infty$. Specifically, by Remark~\ref{rem:explength}, if $0<\eta<1$ and a pairwise uniform inner code yields an increment of at least $\exp(-(2+\eta)q)$ over $B_q$, then its length is at least $\exp((1-\eta)q/2-O(\sqrt q))$. Thus, inner codes of exponential length, as in the proof of Theorem~\ref{thm:exponential}, are necessary for increments on this scale.

\begin{remark*}
An MDS code of dimension two is a code $C\subseteq\cA^\ell$ with $|C|=q^2$ in which any two coordinates determine a codeword. It is pairwise uniform and has length at most $q+1$ (Remark~\ref{rem:limit}), so the bound $B_q+p_q^{3+o(1)}$ applies to it. If ${q\not\equiv2\pmod4}$, the inner code behind Theorem~\ref{thm:asym} is such a code. Theorems~\ref{thm:asym} and~\ref{thm:limit} show that, for ${q\not\equiv2\pmod4}$, the best lower bound on $R_q$ obtainable from Lemma~\ref{lem:concatenation} with MDS inner codes of dimension two is $B_q+p_q^{3+o(1)}$.
\end{remark*}

Finally, for some fixed values of $q$, the extended Reed--Solomon inner code also leads to improved bounds. We derive two exact identities for the number of its separated subsets, one in terms of functions and one in terms of permutations of $\F_q$ (Lemma~\ref{lem:exactcount}). The resulting counts recover the examples for $q=3,4$ in \cite[Remark 4 and Corollary 3.6]{XY2023} and improve the numerical lower bounds reported in \cite[Table 3]{XY2023} for $q=8,9,11$.

The paper is organized as follows. Section~\ref{sec:prelim} collects the concatenation inequality, the Bonferroni inequalities and two counting tools for intersections of separation events. We prove Theorem~\ref{thm:allq} and its refinement Theorem~\ref{thm:length} in Section~\ref{sec:allq}, Theorems~\ref{thm:exponential} and~\ref{thm:general} together with Corollary~\ref{cor:fixedd} in Section~\ref{sec:simplex}, Theorem~\ref{thm:asym} and Corollary~\ref{cor:special} in Section~\ref{sec:asym}, and Theorem~\ref{thm:limit} in Section~\ref{sec:limit}. Section~\ref{sec:counts} contains exact counts for extended Reed--Solomon inner codes, and Section~\ref{sec:remaining} discusses the remaining assertion of the conjecture and further questions.

\section{Preliminaries}\label{sec:prelim}

This section collects the tools used in the later proofs. Section~\ref{sec:concat} recalls the concatenation inequality of Xing and Yuan, converts an upper bound for the failure probability of an inner code into a lower bound for the rate, and states the Bonferroni inequalities. All inner codes in this paper are compared with the independent model, in which the $\ell$ coordinates of the sampled words are independent and uniform. For a pairwise uniform code, intersections of one or two separation events have the same probabilities as in this model, by \eqref{eq:pairwise} below, so the comparison is decided by intersections of three or more events. Section~\ref{sec:inputs} computes the probability of an intersection of three events asymptotically by averaging Eberhard's count of additive triples of bijections, uniformly over all abelian groups $G$ of order $q$ and all subgroups $U\le G$. Section~\ref{sec:entropy} bounds intersections of three or more events from above by Luria's entropy method, uniformly over all 
pairwise uniform codes and all sets of coordinates. This uniformity is essential, since the number of sets of coordinates grows with $q$ in our applications.

Throughout, $N$ denotes the length of a perfect hash code and $\ell$ the length of an inner code. For a prime power $s$, $\F_s$ denotes the field of order $s$, and $\Z/m\Z$ denotes the ring of integers modulo $m$. Logarithms in rate expressions have base two, and $\ln$ denotes the natural logarithm. Entropies are also taken with natural logarithms. Samples from an inner code are always independent and uniform, with replacement. We use Stirling's formula in the form \cite{Robbins1955}
\begin{equation}\label{eq:stirling}
	\sqrt{2\pi n}\,\Bigl(\frac ne\Bigr)^n\le n!
	\le e^{1/(12n)}\sqrt{2\pi n}\,\Bigl(\frac ne\Bigr)^n.
\end{equation}
In particular, $\sqrt{2\pi q}\,e^{-q}\le p_q\le e^{1/(12q)}\sqrt{2\pi q}\,e^{-q}$, and hence $\ln p_q=-q+O(\ln q)$.

\subsection{The concatenation inequality}\label{sec:concat}

\begin{definition}\label{def:separation}
	Let $\cA$ be an alphabet of size $q$, and let $C\subseteq\cA^\ell$. A coordinate $i\in[\ell]$ \emph{separates} words $c_1,\ldots,c_q\in C$ if $\{c_{1,i},\ldots,c_{q,i}\}=\cA$. Let $\cS(C)$ be the family of $q$-element subsets of $C$ that are separated by at least one coordinate. Let $F(C)$ be the probability that no coordinate separates $q$ words sampled from $C$ independently and uniformly, with replacement.
\end{definition}

A separated $q$-tuple contains no repeated word. Hence, the ordered $q$-tuples of words separated by some coordinate are exactly the orderings of the sets in $\cS(C)$, and there are $q!\,|\cS(C)|$ of them. Dividing by the number $|C|^q$ of all $q$-tuples, we obtain
\begin{equation}\label{eq:failure-count}
	F(C)=1-\frac{q!\,|\cS(C)|}{|C|^q}.
\end{equation}
The following concatenation lemma is due to Xing and Yuan. By \eqref{eq:failure-count}, it is a restatement of \cite[Lemma 2.3]{XY2023}, whose proof concatenates $C$ with a randomly chosen outer code.

\begin{lemma}[Xing--Yuan]\label{lem:concatenation}
	If $C\subseteq\cA^\ell$ and $0<F(C)<1$, then
	\[ R_q\ge -\frac{\log_2 F(C)}{\ell(q-1)}.\]
\end{lemma}

Lemma~\ref{lem:concatenation} reduces the construction of good perfect hash codes to the construction of a single inner code. By \eqref{eq:failure-count}, the lower bound in Lemma~\ref{lem:concatenation} equals
\[ -\frac{1}{\ell(q-1)}\log_2\Bigl(1-\frac{q!\,|\cS(C)|}{|C|^q}\Bigr),\]
which increases with the ratio $q!\,|\cS(C)|/|C|^q$. Thus, for given $q$ and $\ell$, we only need to find inner codes $C$ for which $|\cS(C)|$ is as large as possible compared with $|C|^q/q!$. The whole space $\cA^\ell$ serves as the benchmark. Its coordinates are independent, and each of them separates $q$ independent uniform words with probability $p_q$.
Hence, $F(\cA^\ell)=(1-p_q)^\ell$, and for $C=\cA^\ell$ the lower bound in Lemma~\ref{lem:concatenation} is exactly the probabilistic bound $B_q$. Consequently, the lower bound obtained from an inner code $C$ exceeds $B_q$ if and only if $F(C)<(1-p_q)^\ell$, that is, if and only if $q!\,|\cS(C)|/|C|^q>1-(1-p_q)^\ell$. All constructions in this paper are designed to make the difference $(1-p_q)^\ell-F(C)$ positive and as large as possible.

For a code $C\subseteq\cA^\ell$, let
\[ \Phi(C)=-\frac{\log_2F(C)}{\ell(q-1)}.\]
We have $F(C)>0$ for every nonempty code $C$, since two of the $q\ge2$ sampled words coincide with positive probability, and repeated words are not separated. The next lemma converts bounds for $F(C)$ into bounds for $\Phi(C)$.

\begin{lemma}\label{lem:increment}
Let $C\subseteq\cA^\ell$ be a nonempty code over an alphabet of size $q\ge2$.
\begin{enumerate}
\item[(i)] If $F(C)\le(1-p_q)^\ell-X$ for some $0<X<(1-p_q)^\ell$, then $F(C)<1$ and
\[ R_q\ge\Phi(C) \ge B_q+\frac{X}{\ell(q-1)(1-p_q)^\ell\ln2}
\ge B_q+\frac{X}{\ell(q-1)\ln2}.\]
\item[(ii)] If $F(C)\ge(1-p_q)^\ell-Y$ for some $0\le Y<(1-p_q)^\ell$, then
\[ \Phi(C)\le B_q+\frac{Y}{\ell(q-1)\bigl((1-p_q)^\ell-Y\bigr)\ln2}.\]
\end{enumerate}
\end{lemma}

\begin{proof}
By the definition of $B_q$, we have $-\log_2\bigl((1-p_q)^\ell\bigr)=\ell(q-1)B_q$. Let $y=X/(1-p_q)^\ell$ in part (i) and $y=Y/(1-p_q)^\ell$ in part (ii), so that $0\le y<1$ and $-\log_2\bigl((1-p_q)^\ell(1-y)\bigr)=\ell(q-1)B_q-\log_2(1-y)$. In part (i), we have $0<F(C)\le(1-p_q)^\ell(1-y)<1$, so $R_q\ge\Phi(C)$ by Lemma~\ref{lem:concatenation}.
Using the inequality $-\ln(1-y)\ge y$, we obtain
\[ -\log_2F(C)\ge\ell(q-1)B_q-\log_2(1-y) \ge\ell(q-1)B_q+\frac{y}{\ln2}.\]
Dividing by $\ell(q-1)$, we obtain the second inequality in (i). The third follows from $(1-p_q)^\ell\le1$.

In part (ii), we have $F(C)\ge(1-p_q)^\ell(1-y)>0$. Since $-\ln(1-y)\le y/(1-y)$, we have
\[ -\log_2F(C)\le\ell(q-1)B_q-\log_2(1-y) \le\ell(q-1)B_q+\frac{Y}{\bigl((1-p_q)^\ell-Y\bigr)\ln2}.\]
Dividing by $\ell(q-1)$ proves part (ii).
\end{proof}

Suppose that $F(C)=(1-p_q)^\ell-X$ with $X>0$, where $C$ runs through a sequence of codes with $X/(1-p_q)^\ell\to0$. Applying both parts of Lemma~\ref{lem:increment} with $Y=X$, we obtain
\[ \Phi(C)=B_q+(1+o(1))\frac{X}{\ell(q-1)(1-p_q)^\ell\ln2}.\]
The same computation applies to the right-hand sides of \eqref{eq:xy-regular} and \eqref{eq:xy-exceptional}, with $\ell=4$, $X=(3q/\sqrt e+o(q))p_q^3$ and with $\ell=3$, $X=(q/(2\sqrt e)+o(q))p_q^3$, respectively. Since $(1-p_q)^\ell\to1$ for fixed $\ell$, these right-hand sides exceed $B_q$ by $(3/(4\sqrt e\ln2)+o(1))p_q^3$ and $(1/(6\sqrt e\ln2)+o(1))p_q^3$, as stated after \eqref{eq:scale}.

We also use the Bonferroni inequalities \cite{GalambosSimonelli1996}. Let $E_1,\ldots,E_\ell$ be events in a probability space. For $h\ge1$, let $S_h$ be the sum of $\Prb(\bigcap_{j\in J}E_j)$ over all $h$-element subsets $J\subseteq[\ell]$, and let $S_0=1$. Then, for every integer $K\ge0$,
\begin{equation}\label{eq:bonferroni-general}
	(-1)^{K+1}\Bigl(\Prb\Bigl(\bigcap_{j=1}^\ell E_j^{c}\Bigr)
	-\sum_{h=0}^K(-1)^hS_h\Bigr)\ge0.
\end{equation}
In other words, the partial sum $\sum_{h=0}^K(-1)^hS_h$ is an upper bound for the probability that none of the events occurs when $K$ is even, and a lower bound when $K$ is odd. For $K\ge\ell$, the partial sum is the inclusion--exclusion formula and equals this probability.

For $\ell$ independent events of probability $p\in[0,1]$, we have $S_h=\binom\ell hp^h$. Hence, \eqref{eq:bonferroni-general} shows that the truncated binomial sums
\[ P_K(p)=\sum_{h=0}^K(-1)^h\binom\ell hp^h\]
satisfy $(-1)^{K+1}\bigl((1-p)^\ell-P_K(p)\bigr)\ge0$ for every $K\ge0$. Since the difference $P_{K+1}(p)-P_K(p)$ equals $(-1)^{K+1}\binom\ell{K+1}p^{K+1}$, applying this inequality to $K$ and to $K+1$, we obtain
\begin{equation}\label{eq:binomial-truncation}
	0\le(-1)^{K+1}\bigl((1-p)^\ell-P_K(p)\bigr)\le\binom\ell{K+1}p^{K+1}
	\qquad(K\ge0).
\end{equation}

Let $C\subseteq\cA^\ell$, with $\ell\ge2$, be a pairwise uniform code, as defined before Theorem~\ref{thm:limit}. For two distinct coordinates $j$ and $k$, each pair of symbols has exactly $|C|/q^2$ preimages under the projection onto these two coordinates, and hence each symbol has exactly $|C|/q$ preimages under the projection onto either of them. Now sample $q$ words $c_1,\ldots,c_q$ of $C$ independently and uniformly, and let $E_j$ be the event that the $j$th coordinate separates them. The pairs $(c_{i,j},c_{i,k})$, $1\le i\le q$, are independent and uniform on $\cA^2$. Hence, the columns $(c_{1,j},\ldots,c_{q,j})$ and $(c_{1,k},\ldots,c_{q,k})$ are independent, and each of them consists of $q$ independent uniform symbols. Since $E_j$ depends only on the $j$th column, and $q$ independent uniform symbols are all distinct with probability $p_q$, we obtain
\begin{equation}\label{eq:pairwise}
	\Prb(E_j)=p_q~~\text{and}~~\Prb(E_j\cap E_k)=p_q^2\qquad(j\ne k).
\end{equation}

\subsection{Additive triples of bijections}\label{sec:inputs}

In this subsection, $G$ is a finite abelian group of order $q$, and $S_G$ is the set of bijections $[q]\to G$. For $f\colon[q]\to G$, define
\[ N_G(f)=|\{(\pi_1,\pi_2,\pi_3)\in S_G^3:\pi_1+\pi_2+\pi_3=f\}|.\]
The main input for intersections of three separation events is the following theorem of Eberhard.

\begin{theorem}[{\!\!\cite[Theorem 1.3]{Eberhard2017}}]\label{thm:eberhard}
	Let $G$ be a finite abelian group of order $q$, and let $f\colon[q]\to G$ satisfy $\sum_{i\in[q]}f(i)=\sum_{g\in G}g$. Then, as $q\to\infty$,
	\begin{equation}\label{eq:singularseries}
		N_G(f)=\bigl(\mathfrak S_G(f)+o(1)\bigr)\frac{(q!)^3}{q^{q-1}},
		\quad\text{where}\quad
		\mathfrak S_G(f)= \exp\left(-\frac1{2q^2}\sum_{g\in G}|f^{-1}(g)|^2\right).
	\end{equation}
\end{theorem}

We call the functions $f$ in Theorem~\ref{thm:eberhard} \emph{admissible}. Since the group $G$ and the admissible function $f$ in Theorem~\ref{thm:eberhard} may depend on $q$ arbitrarily, the error term in \eqref{eq:singularseries} is uniform over all abelian groups $G$ of order $q$ and all admissible $f$. To see this, let $Q_q=(q!)^3/q^{q-1}$, and let $\alpha(q)$ be the maximum of $|N_G(f)/Q_q-\mathfrak S_G(f)|$ over all pairs $(G,f)$, where $G$ runs through representatives of the finitely many isomorphism classes of abelian groups of order $q$ and $f$ is admissible. Both quantities are invariant under isomorphisms of $G$. Choosing for each $q$ a pair at which the maximum is attained, we obtain a sequence to which Theorem~\ref{thm:eberhard} applies. Hence, $\alpha(q)\to0$, and
\begin{equation}\label{eq:eberhard-uniform}
	\left|\frac{N_G(f)}{Q_q}-\mathfrak S_G(f)\right|\le\alpha(q)
\end{equation}
for every abelian group $G$ of order $q$ and every admissible $f$. This uniformity is essential below, since we average over exponentially many admissible functions.

In particular, if $\sum_{g\in G}g=0$, then $f=0$ is admissible and $\mathfrak S_G(0)=e^{-1/2}$, so that
\begin{equation}\label{eq:zero-triples}
	N_G(0)=\left(e^{-1/2}+o(1)\right)\frac{(q!)^3}{q^{q-1}}.
\end{equation}
The case $f=0$ for odd cyclic groups is due to Eberhard, Manners and Mrazovi\'c \cite{EMM2019}.

The counts $N_G(0)$ are closely related to complete mappings and transversals in group tables. After fixing $\pi_1$ and relabelling its domain by $G$, a solution of $\pi_1+\pi_2+\pi_3=0$ is equivalent to a bijection $\theta\colon G\to G$ for which $x\mapsto x+\theta(x)$ is also bijective. Thus, $N_G(0)/q!$ is the number of complete mappings of $G$. Eberhard, Manners and Mrazovi\'c \cite{EMM2022} obtained an asymptotic count in the broader setting of finite groups satisfying the Hall--Paige condition. For an abelian group $G$, this condition means $\sum_{g\in G}g=0$, and their count recovers \eqref{eq:zero-triples}. M\"uyesser and Pokrovskiy \cite{MP2025} proved a random-set extension of the Hall--Paige theorem. The nonconstant right-hand sides needed below, however, require the more specific formula \eqref{eq:singularseries}.

Let $\sigma_G=\sum_{g\in G}g$. Since $g\mapsto-g$ permutes $G$, we have $\sigma_G=-\sigma_G$, that is, $2\sigma_G=0$. For a subgroup $U\le G$, let $C_U(G)=\{(g_1,g_2,g_3)\in G^3:g_1+g_2+g_3\in U\}$, and let $t_3(G,U)$ be the probability that the three coordinates of $C_U(G)$ separate $q$ independent uniform words of $C_U(G)$. The following lemma averages the singular series in \eqref{eq:singularseries} over all admissible right-hand sides with values in $U$.

\begin{lemma}\label{lem:cosettriples}
Let $G$ be an abelian group of order $q\ge3$, and let $U\le G$ be a subgroup of order $h$. If $\sigma_G\notin U$, then $t_3(G,U)=0$. If $\sigma_G\in U$, then
\[\left|\frac{h\,t_3(G,U)}{q\,p_q^3}-e^{-1/(2h)}\right|\le\alpha(q)+\frac{1-1/h}{2q}.\]
In particular, $t_3(G,U)=(e^{-1/(2h)}+o(1))(q/h)p_q^3$ as $q\to\infty$, uniformly over all such $G$ and $U$.
\end{lemma}

\begin{proof}
	We first express $t_3(G,U)$ through the numbers $N_G(f)$. The code $C_U(G)$ has $hq^2$ words, since $g_1$ and $g_2$ are arbitrary and $g_3$ lies in the coset $U-g_1-g_2$. Let $w_1,\ldots,w_q$ be words of $C_U(G)$, and let $\pi_j\colon[q]\to G$ be the $j$th column, given by $\pi_j(i)=w_{i,j}$. The $q$-tuple is separated by all three coordinates exactly when $\pi_1,\pi_2,\pi_3$ are bijections. The condition that every $w_i$ lies in $C_U(G)$ states that $f=\pi_1+\pi_2+\pi_3$ takes values in $U$. Hence, the number of ordered $q$-tuples of words of $C_U(G)$ separated by all three coordinates is the sum of $N_G(f)$ over all $f\colon[q]\to U$, and
	\[ t_3(G,U)=\frac{1}{(hq^2)^q}\sum_{f\colon[q]\to U}N_G(f).\]
	If $N_G(f)>0$, then $f=\pi_1+\pi_2+\pi_3$ for bijections $\pi_j$, and summing over $[q]$, we obtain $\sum_{i\in[q]}f(i)=3\sigma_G=\sigma_G$. If $\sigma_G\notin U$, no function $f\colon[q]\to U$ has this property, since $\sum_if(i)\in U$. Hence, $t_3(G,U)=0$ in this case.
	
	Assume now that $\sigma_G\in U$, and let $\cF$ be the set of functions $f\colon[q]\to U$ with $\sum_{i\in[q]}f(i)=\sigma_G$. Only the functions in $\cF$ contribute to the sum above, and every $f\in\cF$ is admissible. Since $f(1)$ is determined by $f(2),\ldots,f(q)$, we have $|\cF|=h^{q-1}$. Let $A$ be the average of $\mathfrak S_G(f)$ over $f\in\cF$, and let $\bar e$ be the average over $f\in\cF$ of $N_G(f)/Q_q-\mathfrak S_G(f)$. By \eqref{eq:eberhard-uniform}, we have $|\bar e|\le\alpha(q)$. Since $Q_q=(q!)^3/q^{q-1}$ and $p_q^3=(q!)^3/q^{3q}$, we obtain
	\[ t_3(G,U)=\frac{h^{q-1}}{h^qq^{2q}}\,(A+\bar e)\,Q_q
	=\frac qh\,(A+\bar e)\,p_q^3.\]
	Since $\bar e$ is an average and not a sum, the number $|\cF|$ of terms does not enter its bound.
	
	It remains to estimate $A$. Let $f$ be uniformly distributed on $\cF$, and let $i,j\in[q]$ be distinct. Since $q\ge3$, there is a further position whose value absorbs the sum condition, and hence every choice of $f(i)$ and $f(j)$ has exactly $h^{q-3}$ extensions to a function in $\cF$. Thus, $f(i)$ and $f(j)$ are independent and uniform on $U$. Let $X=\sum_{u\in U}|f^{-1}(u)|^2$, so that $\mathfrak S_G(f)=e^{-X/(2q^2)}$. Since $\sum_{u\in U}|f^{-1}(u)|=q$, we have $X\ge q^2/h$ by the Cauchy--Schwarz inequality. Moreover, for each $u\in U$,
	\[ \E|f^{-1}(u)|^2=\sum_{i,j\in[q]}\Prb\bigl(f(i)=u=f(j)\bigr) =\frac qh+\frac{q(q-1)}{h^2},\]
	and summing over the $h$ elements of $U$, we obtain
	\[ \E X=q+\frac{q(q-1)}h=\frac{q^2}h+q\Bigl(1-\frac1h\Bigr).\]
	Since $X\ge q^2/h$, we have $\mathfrak S_G(f)\le e^{-1/(2h)}$ for every $f\in\cF$, and hence $A\le e^{-1/(2h)}$. For the lower bound, we apply the inequality $|e^{-x/2}-e^{-y/2}|\le|x-y|/2$ for $x,y\ge0$ with $x=X/q^2\ge y=1/h$, and obtain $\mathfrak S_G(f)\ge e^{-1/(2h)}-(X/q^2-1/h)/2$. Averaging over $f\in\cF$ and using the formula for $\E X$, we obtain
	\[ e^{-1/(2h)}-\frac{1-1/h}{2q} \le A\le e^{-1/(2h)}.\]
	Combining this with $|\bar e|\le\alpha(q)$ and the expression for $t_3(G,U)$, we obtain the asserted inequality.
\end{proof}

For $U=\{0\}$, the code $C_U(G)$ is the MDS code $\{(x,y,-x-y):x,y\in G\}$, and Lemma~\ref{lem:cosettriples} reduces to \eqref{eq:zero-triples}. Larger subgroups $U$ occur in Sections~\ref{sec:simplex} and~\ref{sec:asym}, and in Section~\ref{sec:simplex} the order of $U$ may grow with $q$.

\subsection{An entropy bound for separation intersections}\label{sec:entropy}

For intersections of four or more separation events, and for upper bounds on intersections of three, we use Luria's entropy method. Luria \cite[Theorem 3.1]{Luria2017} proved the following. Let $h$ be fixed, and let $H$ be an $h$-uniform, $k$-regular hypergraph on $v$ vertices in which every pair of vertices lies in $o(k)$ common edges. Then, as $k\to\infty$, the number $\operatorname{pm}(H)$ of perfect matchings of $H$ satisfies
\[ \operatorname{pm}(H)\le\left((1+o(1))\frac{k}{e^{h-1}}\right)^{v/h}.\]
We need a version with an explicit error term that applies to all pairwise uniform codes simultaneously. Such a version is not stated in \cite{Luria2017}, and we derive it from Luria's argument, retaining the error terms. The size of the code is not assumed to be a fixed multiple of $q^2$.

\begin{lemma}\label{lem:matching}
Let $C\subseteq\cA^\ell$ be a pairwise uniform code over an alphabet $\cA$ of size $q$, let $J\subseteq[\ell]$ be a set of $h\ge3$ coordinates, and let $t_J(C)$ be the probability that every coordinate in $J$ separates $q$ independent uniform words of $C$. Let
\[ c_h=1+\binom h2~~\text{and}~~\omega_h(q)=\int_0^1\ln\left(1+\frac{c_h}{q t^{h-1}}\right)\,dt.\]
Then we have
\begin{equation}\label{eq:intersectionbound}
t_J(C)\le p_q\exp\bigl(-(h-1)q+q\omega_h(q)\bigr),
\end{equation}
and
\begin{equation}\label{eq:entropy-error-rate}
0\le\omega_h(q)\le\left(h-1+\ln2+\frac{1}{h-2}\right)\left(\frac{c_h}{q}\right)^{1/(h-1)}.
\end{equation}
In particular, for each fixed $h$, we have $t_J(C)\le\exp(-hq+o(q))$ uniformly over all such codes $C$ and sets $J$.
\end{lemma}

\begin{proof}
The proof has three steps. We first express $t_J(C)$ through the number of perfect matchings of a hypergraph. Then we bound this number by Luria's entropy argument, and finally we estimate the integral $\omega_h(q)$.
	
Let $m=|C|$. Since $C$ is pairwise uniform, every pair of symbols in two given coordinates has the same number of preimages, and these $q^2$ numbers sum to $m$. Hence, $\lambda=m/q^2$ is a positive integer, and $k=m/q=q\lambda\ge q$. Let $H$ be the $h$-uniform hypergraph on the vertex set $J\times\cA$ with one edge $e_c=\{(j,c_j):j\in J\}$ for each $c\in C$, where distinct codewords correspond to distinct edges even if their vertex sets coincide. It has $v=hq$ vertices, and every vertex lies in exactly $k$ edges. Two vertices $(j,a)$ and $(j',a')$ with $j\ne j'$ lie in exactly $\lambda$ common edges, and two distinct vertices with the same first component lie in no common edge. A perfect matching of $H$ is a set of $q$ edges covering every vertex exactly once, that is, a $q$-element subset of $C$ separated by every coordinate in $J$. An ordered $q$-tuple of codewords is separated by all coordinates in $J$ exactly when its entries are distinct and form such a subset. Therefore,
\begin{equation}\label{eq:matching-normalization}
t_J(C)=\frac{q!\operatorname{pm}(H)}{m^q}.
\end{equation}
If $\operatorname{pm}(H)=0$, there is nothing to prove, so we assume that $\operatorname{pm}(H)\ge1$.
	
Let $X$ be a uniformly chosen perfect matching of $H$, and for a vertex $w$, let $X_w$ be the edge of $X$ containing $w$. The entropy of $X$ is $\ln\operatorname{pm}(H)$. Assign independent priorities to the vertices, uniform in $[0,1]$ and independent of $X$. Almost surely, the priorities are distinct, and we list the vertices as $w_1,\ldots,w_v$ in decreasing order of priority. Since $X$ determines the family $(X_w)_w$ and is determined by it, it follows from the chain rule for entropy that
\[ \ln\operatorname{pm}(H) =\sum_{i=1}^v H\bigl(X_{w_i}\bigm|X_{w_1},\ldots,X_{w_{i-1}},
\text{priorities}\bigr),\]
where $H(\cdot\mid\cdot)$ denotes conditional entropy. Call $w$ \emph{first} in $X_w$ if $w$ has the largest priority among the $h$ vertices of $X_w$. If $w_i$ is not first in $X_{w_i}$, then $X_{w_i}=X_{w_j}$ for some $j<i$, and the $i$th term vanishes. If $w_i$ is first in $X_{w_i}$, then $X_{w_i}$ is one of the edges through $w_i$ that are disjoint from $X_{w_1},\ldots,X_{w_{i-1}}$. We denote the number of these edges by $N_{w_i}$. Whether $w_i$ is first in $X_{w_i}$, and the number $N_{w_i}$, are determined by the conditioning variables. Since the entropy of a random variable is at most the logarithm of the size of its support, we obtain
\[\ln\operatorname{pm}(H)\le\sum_w\E\bigl[\mathbf 1_{\{w\text{ is first in }X_w\}}\ln N_w\bigr],\]
where the expectation is taken over $X$ and the priorities.
	
Fix $X$. Call an edge \emph{bad} if it contains two vertices of some edge of $X$, and let $b_w(X)$ be the number of bad edges containing $w$. Each of the $q$ edges of $X$ contains $\binom h2$ pairs of vertices, and each such pair lies in at most $\lambda$ edges. Hence, there are at most $q\binom h2\lambda$ bad edges. Since every bad edge contains $h$ vertices,
\begin{equation}\label{eq:bad-incidences}
\sum_w b_w(X)\le hq\binom h2\lambda=v\binom h2\lambda.
\end{equation}
Now fix also a vertex $w$, and condition on the event that the priority of $w$ equals $x$ and that $w$ is first in $X_w$. Given that the priority of $w$ equals $x$, this event requires the other $h-1$ vertices of $X_w$ to have priority less than $x$, and therefore has probability $x^{h-1}$. The edges counted by $N_w$ are $X_w$ itself, some of the bad edges through $w$, and some of the edges through $w$ that are not bad. Let $e$ be an edge through $w$ that is not bad. Then no two vertices of $e$ lie in a common edge of $X$. In particular, the other $h-1$ vertices of $e$ lie in $h-1$ distinct edges of $X$, all different from $X_w$. The edge $e$ is disjoint from the previously exposed edges only if none of these $h-1$ edges of $X$ has been exposed before $w$, that is, only if all their $h(h-1)$ vertices have priority less than $x$. These vertices lie outside $X_w$, so their priorities are independent of the conditioning event, and this happens with probability $x^{h(h-1)}$. Since at most $k$ edges contain $w$, we obtain
\[ \E\bigl(N_w\mid X,\ \text{priority}(w)=x,\ w\text{ first in }X_w\bigr) \le1+b_w(X)+k x^{h(h-1)}.\]
The term $1$ accounts for $X_w$, which is also bad. This double count only weakens the bound. By Jensen's inequality for the concave function $\ln$, the conditional expectation of $\ln N_w$ is at most $\ln(1+b_w(X)+kx^{h(h-1)})$. Integrating over the priority $x$ of $w$ and averaging over $X$, we obtain
\[ \ln\operatorname{pm}(H) \le\E_X\sum_w\int_0^1 x^{h-1}\ln\bigl(1+b_w(X)+kx^{h(h-1)}\bigr)\,dx.\]
For each $c\ge0$, the function $b\mapsto\ln(1+b+c)$ is concave and increasing. Applying Jensen's inequality to the average over the $v$ vertices and then \eqref{eq:bad-incidences}, we may replace each $b_w(X)$ by $\binom h2\lambda$. Hence,
\[ \ln\operatorname{pm}(H) \le v\int_0^1x^{h-1}\ln\left(1+\binom h2\lambda+kx^{h(h-1)}\right)\,dx.\]
Since $k\ge q$ and $\lambda=k/q$, we have $1+\binom h2\lambda\le k(1+\binom h2)/q=kc_h/q$. We substitute $t=x^h$, so that $x^{h-1}\,dx=dt/h$ and $x^{h(h-1)}=t^{h-1}$. Using $v/h=q$ and $\int_0^1\ln t\,dt=-1$, we obtain
\begin{align*}
\ln\operatorname{pm}(H) &\le q\int_0^1\ln\Bigl(kt^{h-1}+\frac{kc_h}q\Bigr)\,dt\\
&=q\int_0^1\Bigl(\ln k+(h-1)\ln t +\ln\Bigl(1+\frac{c_h}{qt^{h-1}}\Bigr)\Bigr)\,dt\\
&=q\bigl(\ln k-(h-1)+\omega_h(q)\bigr).
\end{align*}
Thus, $\operatorname{pm}(H)\le k^qe^{-(h-1)q+q\omega_h(q)}$. Since $k^q/m^q=q^{-q}$, this bound together with \eqref{eq:matching-normalization} implies \eqref{eq:intersectionbound}.
	
To bound $\omega_h(q)$, let $a=(c_h/q)^{1/(h-1)}$, so that $c_h/(qt^{h-1})=(a/t)^{h-1}$, and extend the integration interval to $(0,\infty)$, which can only increase the integral. For $0<t<a$, we have $(a/t)^{h-1}>1$, and hence $\ln\bigl(1+(a/t)^{h-1}\bigr)\le\ln2+(h-1)\ln(a/t)$, whose integral over $(0,a)$ equals $a\ln2+(h-1)a$. For $t\ge a$, the logarithm is at most $(a/t)^{h-1}$, whose integral over $[a,\infty)$ equals $a/(h-2)$ since $h\ge3$. Adding these two integrals proves \eqref{eq:entropy-error-rate}. Finally, for fixed $h$, it follows from \eqref{eq:entropy-error-rate} that $q\omega_h(q)=O_h(q^{1-1/(h-1)})=o(q)$, and $p_q=e^{-q+O(\ln q)}$ by \eqref{eq:stirling}. Together with \eqref{eq:intersectionbound}, this proves the last assertion.
\end{proof}

The error in Lemma~\ref{lem:matching} depends only on $q$ and $h$. Thus, the estimate holds simultaneously for every pairwise uniform code and every set of $h$ coordinates, however many such sets occur, and no union bound over coordinate sets is needed. In Section~\ref{sec:simplex}, we use it only for finitely many fixed values of $h$, after choosing a Bonferroni order. For $h=3$, we have $c_3=4$, and hence, by \eqref{eq:entropy-error-rate},
\begin{equation}\label{eq:omega3}
q\omega_3(q)\le2(3+\ln2)\sqrt q.
\end{equation}

\section{Proof of Theorem~\ref{thm:allq}}\label{sec:allq}

In this section, we prove Theorem~\ref{thm:allq} and its refinement Theorem~\ref{thm:length}. Both follow from Theorem~\ref{thm:spectral}, a non-asymptotic upper bound for the failure probability of a single-parity-check code over a finite-field component of the alphabet.
Since Theorem~\ref{thm:spectral} holds for every $q$ and every admissible choice of its parameters, Theorem~\ref{thm:allq} and the bound \eqref{eq:shortbound} hold for every $q\ge3$, without an unspecified threshold.

The parity condition leaves every proper coordinate projection uniform, so only the intersection of all $\ell$ separation events differs from the independent model. This intersection has probability $p_q^\ell\tau_\ell$, where $\tau_\ell$ is a normalized return probability of a random walk on $\F_s^q$. For odd $\ell$, we bound $\tau_\ell$ from below by Fourier analysis, keeping the contribution of the constant vectors and controlling the negative Fourier coefficients by Parseval's identity. Theorem~\ref{thm:allq} corresponds to the choice $\ell=2q+1$, and Theorem~\ref{thm:length} to the shortest length allowed by the resulting certificate.

\subsection{A non-asymptotic spectral estimate}

Let $s\ge3$ be a prime power dividing $q$, let $r=q/s$, and identify the alphabet with
\[ \cA=\F_s\times[r].\]
Let $\kappa\colon\cA\to\F_s$ be the projection onto the first component. For an odd integer $\ell\ge3$, define
\[ C_{s,\ell}= \left\{(a_1,\ldots,a_\ell)\in\cA^\ell:\sum_{j=1}^{\ell}\kappa(a_j)=0\right\}.\]
Then $|C_{s,\ell}|=q^\ell/s$. If $s=q$, then $C_{s,\ell}$ is a single-parity-check code. For general $s$, the parity condition involves only the $\F_s$-components, and the components in $[r]$ are unrestricted.

\begin{theorem}\label{thm:spectral}
Let $s$, $r$ and $C_{s,\ell}$ be as above, and let
\[ D_{q,s}=\frac{s^q(r!)^s}{q!}\quad\text{and}\quad L_{q,s,\ell}=s-\frac{D_{q,s}-s}{(s-1)^{\ell-2}}.\]
For every odd integer $\ell\ge3$, we have
\begin{equation}\label{eq:failure-upper}
F(C_{s,\ell})\le (1-p_q)^\ell-(L_{q,s,\ell}-1)p_q^\ell.
\end{equation}
If $L_{q,s,\ell}>1$, then the right-hand side of \eqref{eq:failure-upper} is positive and
\[R_q\ge-\frac1{\ell(q-1)}\log_2\left((1-p_q)^\ell -(L_{q,s,\ell}-1)p_q^\ell\right)>B_q.\]
\end{theorem}

\begin{proof}
Sample $q$ words of $C_{s,\ell}$ independently and uniformly, with replacement, and for $j\in[\ell]$, let $E_j$ be the event that the $j$th coordinate separates these words. Let
\[ \Omega=\{v\in\F_s^q:\text{every element of $\F_s$ occurs exactly $r$ times in $v$}\}.\]
A vector in $\Omega$ is determined by the $r$ positions at which each element of $\F_s$ occurs, and hence $|\Omega|=q!/(r!)^s$. Let $V_1,\ldots,V_\ell$ be independent and uniform on $\Omega$, and let $\tau_\ell=s^q\Prb(V_1+\cdots+V_\ell=\bfzero)$. We first show that only the intersection of all $\ell$ events differs from the independent model, and we express this intersection through $\tau_\ell$.
	
\begin{claim}\label{claim1}
For every proper subset $J\subsetneq[\ell]$, we have
\[\Prb\Bigl(\bigcap_{j\in J}E_j\Bigr)=p_q^{|J|}.\]
Moreover, we have
\begin{equation}\label{eq:exactfailure}
F(C_{s,\ell})=(1-p_q)^\ell-(\tau_\ell-1)p_q^\ell.
\end{equation}
\end{claim}
	
\begin{proof}[Proof of Claim~\ref{claim1}]
Let $J\subsetneq[\ell]$, and fix a coordinate $j_0\in[\ell]\setminus J$. Given a word in $\cA^J$, the letters in the coordinates outside $J\cup\{j_0\}$ can be chosen freely. The parity condition then determines the $\F_s$-component of the letter in the coordinate $j_0$, and its $[r]$-component is arbitrary. Hence, every word in $\cA^J$ has exactly $q^{\ell-|J|-1}r=q^{\ell-|J|}/s$ extensions to a word of $C_{s,\ell}$. Since this number does not depend on the word, the projection of a uniform word of $C_{s,\ell}$ onto the coordinates in $J$ is uniform on $\cA^J$. Thus, the entries of the sampled words in these coordinates are independent and uniform on $\cA$, and the first assertion of the claim follows as in the derivation of \eqref{eq:pairwise}.
		
To compute the intersection of all $\ell$ events, let $X$ be a $q\times\ell$ matrix with independent uniform entries in $\cA$. Let $\mathcal R$ be the event that every row of $X$ lies in $C_{s,\ell}$, and let $\mathcal B$ be the event that every column of $X$ is a bijection $[q]\to\cA$. For each row, the $\ell$ elements $\kappa(X_{i,j})$ are independent and uniform on $\F_s$, so their sum is uniform on $\F_s$.
Thus, each row lies in $C_{s,\ell}$ with probability $1/s$, and since the rows are independent, $\Prb(\mathcal R)=s^{-q}$. The columns are also independent, and each of them is a bijection with probability $p_q$, so that $\Prb(\mathcal B)=p_q^\ell$. Conditionally on $\mathcal R$, the rows of $X$ are independent uniform words of $C_{s,\ell}$, since $\mathcal R$ is the intersection of independent events, one for each row. Therefore,
\[ \Prb\Bigl(\bigcap_{j=1}^\ell E_j\Bigr) =\Prb(\mathcal B\mid\mathcal R).\]
A bijection $\pi\colon[q]\to\cA$ with $\kappa\circ\pi=v$ maps, for each $a\in\F_s$, the $r$ positions $i$ with $v_i=a$ bijectively onto $\{a\}\times[r]$. Hence, each vector $v\in\Omega$ is the $\kappa$-image of exactly $(r!)^s$ bijections $[q]\to\cA$. Conditionally on $\mathcal B$, the columns of $X$ are independent uniform bijections, and therefore their $\kappa$-images are independent and uniform on $\Omega$. The event $\mathcal R$ holds exactly when these $\ell$ vectors sum to $\bfzero$. Thus, $\Prb(\mathcal R\mid\mathcal B)=\Prb(V_1+\cdots+V_\ell=\bfzero)$, and by Bayes' formula,
\[\Prb\Bigl(\bigcap_{j=1}^\ell E_j\Bigr)
=\frac{\Prb(\mathcal B)\Prb(\mathcal R\mid\mathcal B)}{\Prb(\mathcal R)}
=s^qp_q^\ell\Prb(V_1+\cdots+V_\ell=\bfzero)=p_q^\ell\tau_\ell.\]
		
By inclusion--exclusion, we have
\[ F(C_{s,\ell})=\sum_{J\subseteq[\ell]}(-1)^{|J|} \Prb\Bigl(\bigcap_{j\in J}E_j\Bigr),\]
where the empty intersection has probability one. By the first assertion and the formula for $\Prb(\bigcap_{j=1}^\ell E_j)$ above, this sum differs from the binomial expansion $(1-p_q)^\ell=\sum_{J\subseteq[\ell]}(-1)^{|J|}p_q^{|J|}$ only in the term indexed by $[\ell]$, which is $(-1)^\ell p_q^\ell\tau_\ell$ instead of $(-1)^\ell p_q^\ell$. Since $\ell$ is odd, the difference is $-(\tau_\ell-1)p_q^\ell$, which proves \eqref{eq:exactfailure}.
\end{proof}
	
It remains to bound the normalized return probability $\tau_\ell$ from below.
	
\begin{claim}\label{claim2}
Let $\psi$ be a nontrivial additive character of $\F_s$, and let $\lambda_u=\E_{V\in\Omega}\psi(u\cdot V)$ for $u\in\F_s^q$. Then every $\lambda_u$ is real and satisfies
\begin{equation}\label{eq:spectralinterval}
-\frac1{s-1}\le\lambda_u\le1.
\end{equation}
In particular, $\lambda_{a\bfone}=1$ for every $a\in\F_s$. Moreover, we have
\begin{equation}\label{eq:momentparseval}
\tau_\ell=\sum_{u\in\F_s^q}\lambda_u^\ell\quad\text{and}\quad
\sum_{u\in\F_s^q}\lambda_u^2 =\frac{s^q}{|\Omega|}=D_{q,s}.
\end{equation}
\end{claim}
	
	\begin{proof}[Proof of Claim~\ref{claim2}]
Let $t\in\F_s^\times$. If $v\in\Omega$, then $tv\in\Omega$, since the positions at which $tv$ equals $a$ are the positions at which $v$ equals $t^{-1}a$. Thus, multiplication by $t$ permutes $\Omega$, and $tV$ is uniform on $\Omega$ whenever $V$ is. Hence, $\lambda_{tu}=\E\psi(u\cdot(tV))=\lambda_u$. For $z\in\F_s$, applying the orthogonality relations for characters \cite[Section~4.1]{TaoVu2006} to the character $t\mapsto\psi(tz)$ of the additive group of $\F_s$, which is trivial exactly when $z=0$, we see that $\sum_{t\in\F_s}\psi(tz)$ equals $s$ if $z=0$ and $0$ otherwise. Consequently, $\sum_{t\in\F_s^\times}\psi(tz)$ equals $s-1$ if $z=0$ and $-1$ otherwise. Averaging the identity $\lambda_u=\lambda_{tu}$ over $t\in\F_s^\times$, we find
\[ \lambda_u =\frac1{s-1}\E_{V\in\Omega}\sum_{t\in\F_s^\times}\psi\bigl(t(u\cdot V)\bigr)
=\frac{s\Prb(u\cdot V=0)-1}{s-1}.\]
Since $0\le\Prb(u\cdot V=0)\le1$, the number $\lambda_u$ is real and satisfies \eqref{eq:spectralinterval}. Since $s>2$, the elements of $\F_s$ sum to zero. Every $V\in\Omega$ contains each element of $\F_s$ exactly $r$ times, and hence $\bfone\cdot V=r\sum_{a\in\F_s}a=0$.
Consequently, we have $(a\bfone)\cdot V=0$ for every $V\in\Omega$, and $\lambda_{a\bfone}=1$ for every $a\in\F_s$.
		
For $w\in\F_s^q$, the map $u\mapsto\psi(u\cdot w)$ is a character of $\F_s^q$, and it is trivial exactly when $w=\bfzero$. Therefore, it follows from the orthogonality relations for characters \cite[Section~4.1]{TaoVu2006} that
\[\frac1{s^q}\sum_{u\in\F_s^q}\psi(u\cdot w)
=\begin{cases}1,&w=\bfzero,\\0,&w\ne\bfzero,\end{cases}
\qquad w\in\F_s^q.\]
Applying this identity with $w=V_1+\cdots+V_\ell$, taking expectations, and using the independence of $V_1,\ldots,V_\ell$, we obtain
\[ \begin{aligned}
\tau_\ell=s^q\Prb(V_1+\cdots+V_\ell=\bfzero)
&=\sum_{u\in\F_s^q}\E\psi\bigl(u\cdot(V_1+\cdots+V_\ell)\bigr)\\
&=\sum_{u\in\F_s^q} \prod_{j=1}^\ell\E\psi(u\cdot V_j)
=\sum_{u\in\F_s^q}\lambda_u^\ell.
\end{aligned}\]
For the second identity, let $V'$ be an independent copy of $V$. Since $\lambda_u$ is real and $\overline{\psi(z)}=\psi(-z)$, we have $\lambda_u^2=\lambda_u\overline{\lambda_u}=\E_{V,V'}\psi(u\cdot(V-V'))$. Applying the same identity with $w=V-V'$, we obtain
\[ \begin{aligned}
\sum_{u\in\F_s^q}\lambda_u^2 &=\E_{V,V'}\sum_{u\in\F_s^q}\psi\bigl(u\cdot(V-V')\bigr)
=s^q\Prb(V=V')=\frac{s^q}{|\Omega|}=D_{q,s}.
\end{aligned}\]
\end{proof}
	
We now bound the negative terms in the odd moment $\tau_\ell$. Since $\ell$ is odd, a term $\lambda_u^\ell$ in the first identity of \eqref{eq:momentparseval} is negative exactly when $\lambda_u<0$, and then it equals $-|\lambda_u|^\ell$. Keeping the $s$ terms $\lambda_{a\bfone}^\ell=1$ and omitting the other nonnegative terms, we obtain
\begin{equation}\label{eq:momentlower}
\begin{aligned}
\tau_\ell&\ge s-\sum_{\lambda_u<0}|\lambda_u|^\ell\ge s-\frac1{(s-1)^{\ell-2}}\sum_{\lambda_u<0}\lambda_u^2\ge s-\frac{D_{q,s}-s}{(s-1)^{\ell-2}}=L_{q,s,\ell}.
\end{aligned}
\end{equation}
The second inequality holds because $|\lambda_u|\le1/(s-1)$ by \eqref{eq:spectralinterval}, and hence $|\lambda_u|^\ell=|\lambda_u|^{\ell-2}\lambda_u^2\le\lambda_u^2/(s-1)^{\ell-2}$, whenever $\lambda_u<0$. The third holds because the vectors $u$ with $\lambda_u<0$ are nonconstant, and by \eqref{eq:momentparseval} the sum of $\lambda_u^2$ over all nonconstant vectors $u$ equals $D_{q,s}-s$. Combining \eqref{eq:momentlower} with \eqref{eq:exactfailure} proves \eqref{eq:failure-upper}. The right-hand side of \eqref{eq:failure-upper} is at least $F(C_{s,\ell})$, which is positive because two sampled words coincide with positive probability. If $L_{q,s,\ell}>1$, then this right-hand side is also strictly smaller than $(1-p_q)^\ell$. Hence, it follows from Lemma~\ref{lem:concatenation} and the monotonicity of $\log_2$ that
\[\begin{aligned}
R_q\ge-\frac{\log_2F(C_{s,\ell})}{\ell(q-1)}
&\ge-\frac1{\ell(q-1)}\log_2\left((1-p_q)^\ell-(L_{q,s,\ell}-1)p_q^\ell\right)\\
&>-\frac{\log_2\bigl((1-p_q)^\ell\bigr)}{\ell(q-1)}=B_q,
\end{aligned}\]
which is the asserted lower bound for $R_q$.
\end{proof}

\begin{proof}[Proof of Theorem~\ref{thm:allq}]
If $q=2$, then any two distinct binary words differ in some coordinate, and this coordinate separates them. Thus, $[2]^N$ is a perfect $2$-hash code, and hence $M(N,2)=2^N$ and $R_2=1$. Since $p_2=1/2$, we also have $B_2=1$.
	
Now let $q\ge3$, and choose a prime-power divisor $s\ge3$ of $q$ as follows. If $q$ has an odd prime divisor, let $s$ be such a divisor. Otherwise, $q$ is a power of two with $q\ge4$, and we take $s=4$. Let $r=q/s$ and $\ell=2q+1$, which is odd. Since $|\Omega|\ge1$, we have $D_{q,s}=s^q/|\Omega|\le s^q$, and therefore
\begin{align*}
L_{q,s,2q+1} &\ge s-\frac{s^q-s}{(s-1)^{2q-1}}\\
&>s-\frac{s^q}{(s-1)^{2q-1}}\\
&=1+(s-1)\left(1-\left(\frac{s}{(s-1)^2}\right)^q\right)\\
&\ge1+2\left(1-\left(\frac34\right)^3\right)=\frac{69}{32}>2.
\end{align*}
The equality in the third line follows from $s^q/(s-1)^{2q-1}=(s-1)\bigl(s/(s-1)^2\bigr)^q$, and the last inequality uses $s-1\ge2$, $q\ge3$ and $s/(s-1)^2\le3/4$ for $s\ge3$. In particular, $L_{q,s,2q+1}-1>1$, and it follows from \eqref{eq:failure-upper} that $F(C_{s,2q+1})\le(1-p_q)^{2q+1}-p_q^{2q+1}$. As in the proof of Theorem~\ref{thm:spectral}, applying Lemma~\ref{lem:concatenation} and the monotonicity of $\log_2$ then yields
\[ \begin{aligned}
R_q &\ge-\frac{\log_2 F(C_{s,2q+1})}{(2q+1)(q-1)}\\
&\ge-\frac1{(2q+1)(q-1)} \log_2\left((1-p_q)^{2q+1}-p_q^{2q+1}\right)\\
&>-\frac{\log_2(1-p_q)}{q-1}=B_q.
\end{aligned}\]
This proves \eqref{eq:allq}. Together with the case $q=2$, it proves \eqref{eq:sharpness}.
\end{proof}

\begin{remark}\label{rem:xy-parity}
Single-parity-check inner codes of arbitrary length already appear in \cite[Corollaries 3.2 and 3.3]{XY2023}. Let $q\ge3$ be a prime power and $s=q$. Then $r=1$, and negating the last coordinate maps $C_{q,\ell}$ onto the $[\ell,\ell-1]$-MDS code $\{(x_1,\ldots,x_{\ell-1},x_1+\cdots+x_{\ell-1})\}$ of \cite[Corollary 3.2]{XY2023} over the additive group of $\F_q$. Since this map relabels the symbols in one coordinate, the two codes have the same separation counts. For odd $\ell$, comparing \eqref{eq:exactfailure} with \cite[Corollaries 3.2 and 3.3]{XY2023} shows that the criterion of \cite[Corollary 3.3]{XY2023} for this code is equivalent to $\tau_\ell>1$. Hence, Theorem~\ref{thm:spectral} verifies this criterion whenever $L_{q,q,\ell}>1$. By the computation in the proof of Theorem~\ref{thm:allq}, which applies to every prime power $s\ge3$ dividing $q$, this holds in particular for $\ell=2q+1$.
\end{remark}

\begin{remark}\label{rem:lengths}
For particular parameters, Theorem~\ref{thm:spectral} can be substantially better than the uniform choice $\ell=2q+1$. For instance, the exact values
\[L_{16,4,5}=\frac{346011088}{212837625}>1\quad\text{and}\quad
L_{18,3,7}=\frac{145637775}{60988928}>1\]
show that the inner lengths $5$ and $7$ already lead to strict improvements for $q=16$ and $q=18$, respectively. More generally, fix a prime power $s\ge3$ and let $q\to\infty$ through multiples of $s$. By Stirling's formula, $(r!)^s\sim(2\pi r)^{s/2}(r/e)^q$ and $q!\sim\sqrt{2\pi q}\,(q/e)^q$. Since $s^q(r/e)^q=(q/e)^q$, it follows that
\[D_{q,s}\sim\frac{(2\pi r)^{s/2}}{\sqrt{2\pi q}}=\frac{(2\pi q)^{(s-1)/2}}{s^{s/2}}.\]
Since $L_{q,s,\ell}>1$ if and only if $(s-1)^{\ell-1}>D_{q,s}-s$, a strict improvement is certified by some odd length $\ell=(s-1)\ln q/(2\ln(s-1))+O_s(1)$. The choice $\ell=2q+1$ in the proof of Theorem~\ref{thm:allq} has the advantage that it requires no optimization and no exceptional cases.
\end{remark}

\begin{remark}[Relation to iterated quasigroups]\label{rem:oldwalk}
The random walk on balanced columns in \cite[Section 4]{Taranenko2018}, which is also related to Taranenko's work on iterated groups \cite[Theorem 6]{Taranenko2021}, leads to a qualitative counterpart of the preceding proof. We describe its specialization to the present setting, together with the additional observation needed for perfect hashing. Let $K$ be the additive subgroup of $\F_s^q$ generated by $\Omega$. The random walk on $K$ with independent increments uniform on $\Omega$ is irreducible, and its stationary distribution is uniform on $K$. It is also aperiodic. If $V\in\Omega$, then $-V\in\Omega$, so the walk can return to $\bfzero$ in two steps. If $a,b\in\F_s^\times$ satisfy $1+a+b=0$, then $V,aV,bV\in\Omega$ and $V+aV+bV=\bfzero$, so the walk can also return in three steps. Such $a,b$ exist because $s>2$. Since every vector in $\Omega$ has coordinate sum zero, we have $K\subseteq\{v\in\F_s^q:\sum_iv_i=0\}$ and $|K|\le s^{q-1}$. By the convergence theorem for finite Markov chains \cite{LPW2009}, we obtain $\tau_\ell\to s^q/|K|\ge s>1$ as $\ell\to\infty$. Together with \eqref{eq:exactfailure} and Lemma~\ref{lem:concatenation}, this already implies $R_q>B_q$ for every $q\ge3$, but without an explicit inner length. In contrast, the inequality \eqref{eq:momentlower} and the resulting length bounds do not depend on an unspecified mixing threshold.
\end{remark}

\subsection{Shorter effective inner lengths}\label{sec:length}

From the spectral certificate, we also obtain a shorter computable inner length. Here $D_{q,s}$ and $L_{q,s,\ell}$ are as in Theorem~\ref{thm:spectral}. The proof bounds $D_{q,s}$ from above in two ways, by $e^q$ and by a power of $q+1$, and then determines an odd length $\ell$ for which the negative part of \eqref{eq:momentlower} is smaller than $1/(s-1)$.

\begin{theorem}\label{thm:length}
For every integer $q\ge3$, let $\ell_*(q)$ be the smallest odd integer $\ell\ge3$ for which $L_{q,s,\ell}\ge2$ for at least one prime-power divisor $s\ge3$ of $q$. This length is finite, can be computed by exact integer arithmetic, and satisfies
\begin{equation}\label{eq:shortbound}
R_q\ge-\frac1{\ell_*(q)(q-1)}
\log_2\left((1-p_q)^{\ell_*(q)}-p_q^{\ell_*(q)}\right)>B_q.
\end{equation}
Moreover, $\ell_*(q)\le(1+o(1))q/\ln q$ as $q\to\infty$ through all integers. For every prime-power divisor $s\ge3$ of $q$, an explicit upper bound is
\begin{equation}\label{eq:closedlength}
\ell_*(q)\le 2\left\lceil
\frac{\min\{q,(s-1)\ln(q+1)\}}{2\ln(s-1)}
\right\rceil+3.
\end{equation}
\end{theorem}

\begin{proof}
We first bound $D_{q,s}$. By the multinomial theorem, $s^q$ is the sum of the multinomial coefficients $q!/(k_1!\cdots k_s!)$ over the at most $(q+1)^{s-1}$ compositions $k_1+\cdots+k_s=q$, and the balanced coefficient $q!/(r!)^s$ is the largest of them, since moving one unit from a part $k_i$ to a part $k_j\le k_i-2$ multiplies the coefficient by $k_i/(k_j+1)>1$, and repeated moves of this kind lead to $(r,\ldots,r)$ because $s$ divides $q$. Hence, $s^q\le(q+1)^{s-1}q!/(r!)^s$, that is, $D_{q,s}\le(q+1)^{s-1}$. Moreover, $(r!)^s\le r^{rs}=r^q$, and $q!\ge(q/e)^q$ by \eqref{eq:stirling}, so that $D_{q,s}\le s^qr^q(e/q)^q=e^q$. Thus, $D_{q,s}\le e^M$, where $M=\min\{q,(s-1)\ln(q+1)\}$.
	
Let $\ell$ be the odd integer on the right of \eqref{eq:closedlength}. Since $2\lceil x\rceil\ge2x$, we have $\ell\ge3+M/\ln(s-1)$, and hence $(s-1)^{\ell-2}\ge(s-1)^{1+M/\ln(s-1)}=(s-1)e^M$. Thus,
\[ \frac{D_{q,s}-s}{(s-1)^{\ell-2}} <\frac{D_{q,s}}{(s-1)e^M}\le\frac1{s-1}\quad\text{and}\quad L_{q,s,\ell}>s-\frac1{s-1}\ge\frac52>2,\]
where the last step uses $s\ge3$. This proves that $\ell_*(q)$ is finite and satisfies \eqref{eq:closedlength}. Multiplying by $(s-1)^{\ell-2}q!$, we see that the condition $L_{q,s,\ell}\ge2$ is equivalent to $s^q(r!)^s-sq!\le(s-2)(s-1)^{\ell-2}q!$, so it can be tested in exact integer arithmetic. For $\ell=\ell_*(q)$ and a divisor $s$ with $L_{q,s,\ell}\ge2$, we have $F(C_{s,\ell})\le(1-p_q)^\ell-p_q^\ell$ by \eqref{eq:failure-upper}. As in the proof of Theorem~\ref{thm:spectral}, we obtain \eqref{eq:shortbound} from Lemma~\ref{lem:concatenation} and the monotonicity of $\log_2$.
	
For the asymptotic bound on $\ell_*(q)$, fix any prime-power divisor $s\ge3$ of $q$. Such a divisor exists by the proof of Theorem~\ref{thm:allq}. The right-hand side of \eqref{eq:closedlength} is at most $M/\ln(s-1)+5$. We consider two cases. If $s\le q/(\ln q)^3$, then $M\le(s-1)\ln(q+1)$ and $\ln(s-1)\ge\ln2$, so the right-hand side of \eqref{eq:closedlength} is at most
\[ 5+\frac{(s-1)\ln(q+1)}{\ln2} \le5+\frac{q\ln(q+1)}{(\ln q)^3\ln2} =O\left(\frac q{(\ln q)^2}\right).\]
If $s>q/(\ln q)^3$, then $M\le q$ and $\ln(s-1)\ge\ln\bigl(q/(\ln q)^3-1\bigr)=\ln q-3\ln\ln q+o(1)$, so the right-hand side of \eqref{eq:closedlength} is at most
\[ 5+\frac q{\ln(s-1)} \le5+\frac q{\ln q-3\ln\ln q+o(1)} =(1+o(1))\frac q{\ln q}.\]
In both cases, $\ell_*(q)\le(1+o(1))q/\ln q$, and the estimates are uniform in the choice of $s$.
\end{proof}

\begin{remark}
Along primes $q$, the only admissible divisor is $s=q$, and $D_{q,q}=q^q/q!$. Taking logarithms in $L_{q,q,\ell}\ge2$ and using Stirling's formula, we see that $\ell_*(q)\sim q/\ln q$ along primes. This identifies the worst-case order of this particular sufficient spectral certificate, but not a lower bound on the length of the shortest improving inner code, since inner codes of length three or four already lead to strict improvements for all sufficiently large $q$, by \eqref{eq:xy-regular} and \eqref{eq:xy-exceptional}.
\end{remark}

Table~\ref{tab:lengths} compares $\ell_*(q)$ with the length $2q+1$
in Theorem~\ref{thm:allq}. For each $q$, all prime-power divisors
$s\ge3$ are considered, and the displayed $s$ is the smallest one
satisfying $L_{q,s,\ell_*(q)}\ge2$.

\begin{table}[t]
\centering
\setlength{\belowcaptionskip}{3pt}
\caption{Certified inner lengths.}
\label{tab:lengths}
\renewcommand{\arraystretch}{1.25}
\begin{tabular*}{0.6\textwidth}{@{\extracolsep{\fill}}cccc@{}}
\toprule
$q$ & $s$ & $\ell_*(q)$ & $2q+1$\\
\midrule
16  & 4   & 7  & 33\\
18  & 3   & 7  & 37\\
100 & 4   & 9  & 201\\
101 & 101 & 23 & 203\\
210 & 5   & 9  & 421\\
499 & 499 & 81 & 999\\
500 & 4   & 11 & 1001\\
\bottomrule
\end{tabular*}
\end{table}

\section{Proofs of Theorems~\ref{thm:exponential} and~\ref{thm:general}}\label{sec:simplex}

In this section, we prove Theorem~\ref{thm:general} and deduce Theorem~\ref{thm:exponential} and Corollary~\ref{cor:fixedd} from it. The inner code is $C(q,s,d)$, whose coordinates are indexed by the projective points in $\cP_{s,d}$. We show that $F(C(q,s,d))$ is smaller than $(1-p_q)^\ell$ by at least $(1-o(1))Ta_{q,s}p_q^3$, where $T$ is the number of collinear triples of points, and then apply Lemma~\ref{lem:increment}.

In Section~\ref{sec:simplex-triples}, we show that linearly independent sets of coordinates behave exactly as in the independent model, and that every collinear triple has the same separation probability, which we evaluate by Lemma~\ref{lem:cosettriples}. In Section~\ref{sec:simplex-bonferroni}, we truncate inclusion--exclusion at a fixed even order $K$ with $(K-1)\varepsilon>1$. The remaining terms of orders four to $K$ come from linearly dependent sets of coordinates. We count these sets projectively and bound their contributions with the help of Lemma~\ref{lem:matching}. Under \eqref{eq:lengthcondition}, these terms and the truncation error are negligible compared with $Ta_{q,s}p_q^3$. Section~\ref{sec:exp} contains the deduction of Theorem~\ref{thm:exponential} and Corollary~\ref{cor:fixedd}, and Remark~\ref{rem:rs} identifies the extended Reed--Solomon code as a special case.

We use the following terminology. A projective point is a one-dimensional subspace of $\F_s^d$, and we identify it with its representative in $\cP_{s,d}$. The rank of a set of projective points is the dimension of the span of their representatives. A set of $h$ points is independent if its rank is $h$ and dependent otherwise. Three distinct points are collinear exactly when their rank is two.

\subsection{Coordinate projections and triple intersections}\label{sec:simplex-triples}

Let $E_v$ be the event that the coordinate $v\in\cP_{s,d}$ separates $q$ independent uniform words of $C(q,s,d)$. Let $J\subseteq\cP_{s,d}$ have cardinality $h$ and rank $\rho$. The projection of $C(q,s,d)$ onto the coordinates in $J$ maps the word with parameters $u\in\F_s^d$ and $(b_v)_{v\in\cP_{s,d}}$ to $\bigl((u\cdot v,b_v)\bigr)_{v\in J}$. Its image is
\[C_J=\bigl\{\bigl((u\cdot v,b_v)\bigr)_{v\in J}:u\in\F_s^d,\ (b_v)_{v\in J}\in H^J\bigr\}.\]
The linear map $\F_s^d\to\F_s^J$, $u\mapsto(u\cdot v)_{v\in J}$, has rank $\rho$, since its rank equals the dimension of the span of the vectors $v\in J$. Hence, its image has $s^\rho$ elements, and each of them has $s^{d-\rho}$ preimages. The components $b_v$ with $v\in J$ are arbitrary, and the $\ell-h$ components $b_v$ with $v\notin J$ do not affect the projection. Therefore, $|C_J|=s^\rho r^h$, and every word of $C_J$ has exactly $s^{d-\rho}r^{\ell-h}$ preimages. In particular, the projection of a uniform word of $C(q,s,d)$ onto $J$ is uniform on $C_J$. If $\rho=h$, then $C_J=(\F_s\times H)^J$, so the entries of the sampled words in the coordinates of $J$ are independent and uniform. Since each $E_v$ depends only on the $v$th column, it follows that
\begin{equation}\label{eq:independent}
\Prb\Bigl(\bigcap_{v\in J}E_v\Bigr)=p_q^h \quad\text{whenever $J$ is independent.}
\end{equation}
Any two distinct projective points are linearly independent. Taking $|J|=\rho=2$ above, we see that $C(q,s,d)$ is pairwise uniform.

For three distinct points, the rank is either two or three, and rank two means that the points are collinear. Let $v_1,v_2,v_3$ be collinear. Then there exist $c_1,c_2,c_3\in\F_s$ (not all zero) such that $c_1v_1+c_2v_2+c_3v_3=0$. All three coefficients are nonzero, since any two distinct projective points are linearly independent. In the $i$th of these three coordinates, the map $(x,b)\mapsto(c_ix,b)$ is a bijection of $\F_s\times H$, so this coordinatewise relabelling preserves separation.
After relabelling, the projected image consists of the triples $((c_1u\cdot v_1,b_1),(c_2u\cdot v_2,b_2),(c_3u\cdot v_3,b_3))$. Their $\F_s$-components sum to $u\cdot(c_1v_1+c_2v_2+c_3v_3)=0$, and by the preceding paragraph the image has $s^2r^3$ elements. The set of all triples in $(\F_s\times H)^3$ whose $\F_s$-components sum to zero also has $s^2r^3$ elements. Hence, the relabelled image is this set, which is the code $C_U(G)$ of Lemma~\ref{lem:cosettriples} with
\[G=\F_s\times H,~ U=\{0\}\times H~\text{and}~|U|=r.\]
Since the projected distribution is uniform, the probability that the three coordinates $v_1,v_2,v_3$ separate $q$ independent uniform words of $C(q,s,d)$ equals $t_3(G,U)$. In particular, it does not depend on the collinear triple, and we denote it by $t_3(q,s)$.

\begin{lemma}\label{lem:triple}
Let $s\ge3$ be a prime-power divisor of $q$, let $r=q/s$, and write $a_{q,s}=se^{-1/(2r)}-1$. Then $t_3(q,s)=\bigl(se^{-1/(2r)}+o(s)\bigr)p_q^3$ and
\begin{equation}\label{eq:delta}
t_3(q,s)-p_q^3=(1+o(1))a_{q,s}p_q^3\quad\text{and}\quad\Bigl(e^{-1/2}-\frac13\Bigr)s\le a_{q,s}\le s.
\end{equation}
More precisely, let $\xi_q=\alpha(q)+1/(2q)$, where $\alpha(q)$ is as in \eqref{eq:eberhard-uniform}. Then the terms $o(s)$ and $o(1)$ above are bounded in absolute value by $\xi_q s$ and $\xi_q/(e^{-1/2}-1/3)$, respectively. Since $\xi_q$ depends only on $q$ and tends to zero, both estimates are uniform over all prime-power divisors $s\ge3$ of $q$.
\end{lemma}

\begin{proof}
Since $s>2$, the elements of $\F_s$ sum to zero. Hence, we have
\[ \sigma_G=\sum_{(x,b)\in\F_s\times H}(x,b) =\Bigl(r\sum_{x\in\F_s}x,\ s\sum_{b\in H}b\Bigr)
=\Bigl(0,\ s\sum_{b\in H}b\Bigr)\in U.\]
Thus, Lemma~\ref{lem:cosettriples} applies with $h=|U|=r$ and $q/h=s$, and we obtain
\[ \left|\frac{t_3(q,s)}{sp_q^3}-e^{-1/(2r)}\right| \le\alpha(q)+\frac{1-1/r}{2q}\le\xi_q.\]
Multiplying by $sp_q^3$, we see that the term $o(s)$ is bounded in absolute value by $\xi_qs$, uniformly in $s$ and $r$. Furthermore, since $r\ge1$ and $s\ge3$, we have $s\ge a_{q,s}\ge se^{-1/2}-1 \ge\left(e^{-1/2}-\frac13\right)s.$
Thus, $s\xi_q\le\xi_q a_{q,s}/(e^{-1/2}-1/3)$, and subtracting $p_q^3$ from $t_3(q,s)$ proves \eqref{eq:delta} with the stated bound on the error term.
\end{proof}

Every pair of distinct points of $\cP_{s,d}$ spans a two-dimensional subspace, and hence lies on a unique projective line, which has exactly $(s^2-1)/(s-1)=s+1$ points. Choosing two points and then a third point on their line counts each collinear triple three times, once for each of its pairs. Therefore, the number of collinear triples is 
\begin{equation}\label{eq:T3}
T=\frac13\binom\ell2(s-1)  =\frac{\ell(\ell-1)(s-1)}6.
\end{equation}
For $h\ge1$, let $S_h$ denote the sum of the probabilities $\Prb(\bigcap_{v\in J}E_v)$ over all $h$-subsets $J$ of $\cP_{s,d}$, and let $S_0=1$. Then $S_h=0$ for $h>\ell$. Every set of one or two points is independent, and a set of three points is either independent or collinear. By \eqref{eq:independent} and the definition of $t_3(q,s)$,
\[S_1=\ell p_q,\quad S_2=\binom\ell2p_q^2\quad\text{and}\quad S_3=\binom\ell3p_q^3+T\bigl(t_3(q,s)-p_q^3\bigr).\]

\subsection{Dependent sets and the Bonferroni remainder}\label{sec:simplex-bonferroni}

Throughout this subsection, fix $\varepsilon\in(0,1)$, let $s$, $r$ and $d$ be as in Theorem~\ref{thm:general}, and write $C=C(q,s,d)$, $\cP=\cP_{s,d}$ and $p=p_q$. All asymptotic statements refer to $q\to\infty$, and they are uniform over all prime-power divisors $s\ge3$ of $q$ and all $d\ge2$ satisfying \eqref{eq:lengthcondition}. Let
\[\Delta=T\bigl(t_3(q,s)-p^3\bigr).\]

We first determine the size of $\Delta$. Since $\ell\ge s+1\ge4$ and $s\ge3$, we have $\ell-1\ge3\ell/4$ and $s-1\ge2s/3$, and hence $\ell^2s/12\le T\le\ell^2s/6$ by \eqref{eq:T3}. By Lemma~\ref{lem:triple}, we have $\Delta=(1+o(1))Ta_{q,s}p^3$ and $(e^{-1/2}-1/3)s\le a_{q,s}\le s$. Since the factor $1+o(1)$ lies between $1/2$ and $6$ for all sufficiently large $q$, it follows that, for all such $q$,
\begin{equation}\label{eq:Deltascale}
\Delta=(1+o(1))Ta_{q,s}p^3\quad\text{and}\quad c_0\ell^2s^2p^3\le\Delta\le\ell^2s^2p^3,
\end{equation}
where $c_0=(e^{-1/2}-1/3)/24>0$. In particular, $\Delta>0$. Since $3\le s\le q$ and $\ln p=-q+O(\ln q)$ by \eqref{eq:stirling}, we also obtain $\ln\Delta=2\ln\ell-3q+O(\ln q)$.

\begin{claim}\label{claim3}
For $h\ge1$, let $D_h$ be the number of dependent $h$-subsets of $\cP$, and let $U_h=S_h-\binom\ell hp^h$. Then
\begin{equation}\label{eq:dependent}
D_h\le\binom\ell{h-1}\frac{s^{h-1}-1}{s-1} \le2\ell^{h-1}s^{h-2}
\end{equation}
and
\begin{equation}\label{eq:dependent-exact}
U_h=\sum_{\substack{J\subseteq\cP,\ |J|=h\\ J\text{ dependent}}}
\Bigl(\Prb\Bigl(\bigcap_{v\in J}E_v\Bigr)-p^h\Bigr).
\end{equation}
In particular, $U_1=U_2=0$ and $U_3=\Delta$.
\end{claim}

\begin{proof}[Proof of Claim~\ref{claim3}]
Let $J$ be a dependent $h$-subset of $\cP$. Since $J$ is dependent, some $v_0\in J$ has a nonzero coefficient in a nontrivial linear relation among the points of $J$, and hence $v_0$ lies in the span of $J\setminus\{v_0\}$. This span has dimension at most $h-1$, so it contains at most $(s^{h-1}-1)/(s-1)$ projective points. Hence, every dependent $h$-subset arises by choosing $h-1$ points of $\cP$ and then one further point in their span, and the first inequality in \eqref{eq:dependent} follows. Some dependent sets may be counted more than once, which only weakens the bound. The second inequality follows from $\binom\ell{h-1}\le\ell^{h-1}$ and $(s^{h-1}-1)/(s-1)\le s^{h-2}\cdot s/(s-1)\le3s^{h-2}/2$. Both inequalities also hold when $h>d$ or $h>\ell$.
	
The number of $h$-subsets of $\cP$ is $\binom\ell h$, and by \eqref{eq:independent} every independent $h$-subset contributes exactly $p^h$ to $S_h$. Subtracting $\binom\ell hp^h$ from $S_h$ therefore leaves only the dependent $h$-subsets, which proves \eqref{eq:dependent-exact}.
Every set of one or two points is independent, so $U_1=U_2=0$. The dependent triples are the $T$ collinear triples, each of which contributes $t_3(q,s)-p^3$, so $U_3=\Delta$.
\end{proof}

From now on, let $K\ge4$ be an even integer such that $(K-1)\varepsilon>1$, and keep $K$ fixed as $q\to\infty$. For $4\le h\le K$, let
\[ E_h= \begin{cases}
	D_h p\exp\bigl(-(h-1)q+q\omega_h(q)\bigr),&h\text{ even},\\
	D_h p^h,&h\text{ odd}.
\end{cases}\]

\begin{claim}\label{claim4}
For every $h$ with $4\le h\le K$, we have $(-1)^hU_h\le E_h$ and
\[ E_h\le\Delta\exp\bigl(-(h-3)\varepsilon q+o_K(q)\bigr).\]
In particular, $\sum_{h=4}^KE_h=o(\Delta)$.
\end{claim}

\begin{proof}[Proof of Claim~\ref{claim4}]
For even $h$, we drop the nonpositive terms $-p^h$ in \eqref{eq:dependent-exact}. Applying Lemma~\ref{lem:matching} to the pairwise uniform code $C$ and each dependent $h$-subset $J$, we obtain
\[ (-1)^hU_h=U_h \le\sum_{\substack{|J|=h\\ J\text{ dependent}}}
\Prb\Bigl(\bigcap_{v\in J}E_v\Bigr)
\le D_h p\exp\bigl(-(h-1)q+q\omega_h(q)\bigr)=E_h.\]
For odd $h$, since the intersection probabilities are nonnegative, it follows from \eqref{eq:dependent-exact} that $(-1)^hU_h=-U_h\le D_hp^h=E_h$.
If $D_h=0$, then $E_h=0$, and the asserted bound for $E_h$ holds trivially. Let $D_h>0$. By \eqref{eq:dependent}, $s\le q$ and $\ln p=-q+O(\ln q)$, we have $\ln E_h\le(h-1)\ln\ell-hq+O_K(\ln q)+q\omega_h(q)$ for even $h$, and the same bound without the term $q\omega_h(q)$ for odd $h$. By \eqref{eq:entropy-error-rate},
\[ \max_{4\le h\le K}q\omega_h(q)   =O_K\bigl(q^{1-1/(K-1)}\bigr)=o_K(q).\]
Subtracting $\ln\Delta=2\ln\ell-3q+O(\ln q)$ and using $\ln\ell\le(1-\varepsilon)q$, we obtain
\[ \ln\frac{E_h}{\Delta}\le(h-3)(\ln\ell-q)+o_K(q)\le-(h-3)\varepsilon q+o_K(q),\]
which is the asserted bound for $E_h$. Since $h\ge4$ and $K$ is fixed, the sum of the finitely many terms $E_h$ is $o(\Delta)$, uniformly in $s$ and $d$.
\end{proof}

\begin{claim}\label{claim5}
Let $P_K(p)=\sum_{h=0}^{K}(-1)^h\binom\ell hp^h$ and $E_0=P_K(p)-(1-p)^\ell$. Then
\[ 0\le E_0\le\binom\ell{K+1}p^{K+1}\]
and $E_0=o(\Delta)$.
\end{claim}

\begin{proof}[Proof of Claim~\ref{claim5}]
Since $K$ is even, the two inequalities for $E_0$ follow from \eqref{eq:binomial-truncation}. If $\ell\le K$, then $E_0=0$. Otherwise, using the upper bound for $E_0$, the inequality $\binom\ell{K+1}\le\ell^{K+1}$ and \eqref{eq:Deltascale}, we obtain
\[\frac{E_0}{\Delta}\le\frac{\ell^{K+1}p^{K+1}}{c_0\ell^2s^2p^3}
\le\frac{\ell^{K-1}p^{K-2}}{c_0}.\]
By \eqref{eq:lengthcondition} and $\ln p=-q+O(\ln q)$, the logarithm of the right-hand side is at most
\[(K-1)(1-\varepsilon)q-(K-2)q+O_K(\ln q)=\bigl(1-(K-1)\varepsilon\bigr)q+O_K(\ln q).\]
Since $(K-1)\varepsilon>1$, the coefficient of $q$ is negative, and hence $E_0=o(\Delta)$.
\end{proof}

\begin{proof}[Proof of Theorem~\ref{thm:general}]
Since $K$ is even, it follows from the Bonferroni inequality \eqref{eq:bonferroni-general} that $F(C)\le\sum_{h=0}^{K}(-1)^hS_h$. Writing $S_h=\binom\ell hp^h+U_h$ for $1\le h\le K$ and using $U_1=U_2=0$ and $U_3=\Delta$ from Claim~\ref{claim3}, we obtain the exact identity
\[ \sum_{h=0}^{K}(-1)^hS_h  =(1-p)^\ell+E_0-\Delta+\sum_{h=4}^{K}(-1)^hU_h.\]
By Claim~\ref{claim4}, each term $(-1)^hU_h$ with $4\le h\le K$ on the right-hand side is at most $E_h$. Hence, $F(C)\le(1-p)^\ell-X$, where $X=\Delta-E_0-\sum_{h=4}^{K}E_h$.
By Claims~\ref{claim4} and~\ref{claim5}, we have $X=(1-o(1))\Delta$, and hence $X=(1-o(1))Ta_{q,s}p^3$ by \eqref{eq:Deltascale}. In particular, $X>0$ for all sufficiently large $q$. To compare $X$ with $(1-p)^\ell$, note that, by \eqref{eq:lengthcondition} and \eqref{eq:Deltascale},
\[ \ell p\le\exp(-\varepsilon q+O(\ln q))\quad\text{and}\quad
X\le\Delta\le\ell^2s^2p^3\le\exp(-(1+2\varepsilon)q+O(\ln q)).\]
Since $(1-p)^\ell\ge1-\ell p$ by Bernoulli's inequality, we obtain $0<X<(1-p)^\ell$ for all sufficiently large $q$. Thus, part (i) of Lemma~\ref{lem:increment} applies to the bound $F(C)\le(1-p)^\ell-X$. Combining it with the asymptotic formula for $X$, we obtain
\[R_q\ge B_q+\frac{X}{\ell(q-1)\ln2} =B_q+(1-o(1))\frac{Ta_{q,s}}{\ell(q-1)\ln2}p^3.\]
Substituting $T/\ell=(\ell-1)(s-1)/6$ from \eqref{eq:T3} proves \eqref{eq:generalrate}. All error estimates depend only on $q$ and on $\varepsilon$, through the fixed integer $K$. This establishes the stated uniformity.
\end{proof}

\begin{remark}\label{rem:orderfour}
The Bonferroni order $K=4$ suffices when $\ln\ell\le(2/3-\delta)q$ for a fixed $\delta>0$. Its remainder in the independent model is at most $\binom\ell5p_q^5$, and not $\binom\ell4p_q^4$, because the independent four-coordinate intersections cancel exactly. Larger fixed even orders allow every fixed exponent $1-\varepsilon<1$ in \eqref{eq:lengthcondition}. No estimate in this proof is asserted to be uniform when $K$ grows with $q$.
\end{remark}

\subsection{Proof of Theorem~\ref{thm:exponential} and Corollary~\ref{cor:fixedd}}\label{sec:exp}

Theorem~\ref{thm:general} applies to every dimension $d$ with $\ell\le e^{(1-\varepsilon)q}$, and the increment in \eqref{eq:generalrate} grows linearly in $\ell$. For Theorem~\ref{thm:exponential}, we therefore take $d$ as large as this condition allows. Then $\ell$ is close to $e^{(1-\varepsilon)q}$, and the increment is of order $e^{(1-\varepsilon)q}p_q^3$, which is $e^{-(2+\varepsilon)q}$ up to a power of $q$. For Corollary~\ref{cor:fixedd}, we keep $d$ fixed and evaluate the coefficient in \eqref{eq:generalrate} asymptotically.

\begin{proof}[Proof of Theorem~\ref{thm:exponential}]
If $\eta\ge1$, then $e^{-(2+\eta)q}\le e^{-(2+1/2)q}$, so it suffices to prove \eqref{eq:exponential} for $0<\eta<1$. Fix such an $\eta$, let $\varepsilon=\eta/2$, and choose a prime-power divisor $s\ge3$ of $q$ as in the proof of Theorem~\ref{thm:allq}. For $d\ge1$, let $\ell_d=(s^d-1)/(s-1)$. Since $\ell_2=s+1\le q+1\le e^{(1-\varepsilon)q}$ for all sufficiently large $q$, there is a largest integer $d\ge2$ with $\ell_d\le e^{(1-\varepsilon)q}$. We use this $d$, so that \eqref{eq:lengthcondition} holds with $\ell=\ell_d$, and $s\ell+1=\ell_{d+1}>e^{(1-\varepsilon)q}$ by the maximality of $d$.

Since $\ell-1\ge\ell/2$, $s-1\ge2s/3$ and $a_{q,s}\ge(e^{-1/2}-1/3)s\ge s/4$ by \eqref{eq:delta}, the coefficient $(\ell-1)(s-1)a_{q,s}/(6(q-1)\ln2)$ in \eqref{eq:generalrate} is at least $\ell s^2/(72q\ln2)$. Since $s\ell>e^{(1-\varepsilon)q}-1$ and $s\ge3$, this coefficient is at least $(e^{(1-\varepsilon)q}-1)/(24q\ln2)=e^{(1-\varepsilon)q-O(\ln q)}$. By Theorem~\ref{thm:general}, applied with this value of $\varepsilon$, the factor $1-o(1)$ in \eqref{eq:generalrate} depends only on $\varepsilon$ and $q$, so it is at least $1/2$ for all sufficiently large $q$, whatever $s$ and $d$ are. Since also $p_q^3=e^{-3q+O(\ln q)}$ by \eqref{eq:stirling}, we obtain $R_q-B_q\ge e^{-(2+\varepsilon)q-O(\ln q)}$ for all sufficiently large $q$. Since $\varepsilon<\eta$, the right-hand side exceeds $e^{-(2+\eta)q}$ for all sufficiently large $q$. This proves \eqref{eq:exponential}.

By \eqref{eq:exponential}, the lower limit in \eqref{eq:liminf} is at least $-(2+\eta)$ for every $\eta>0$, and letting $\eta$ decrease to zero proves \eqref{eq:liminf}. Finally, $q^Ap_q^3=e^{-3q+O_A(\ln q)}$ by \eqref{eq:stirling}. Hence, by \eqref{eq:exponential} with $\eta=1/2$, we have $(R_q-B_q)/(q^Ap_q^3)\ge e^{q/2-O_A(\ln q)}$, which tends to infinity for every fixed $A>0$. This is \eqref{eq:superpolynomial}.
\end{proof}

\begin{remark}
The proof fixes $\eta$, and hence the Bonferroni order $K$, before letting $q\to\infty$. It does not substitute a function $\varepsilon(q)\to0$ into Theorem~\ref{thm:general}, and it does not use an intersection estimate that is uniform in an unbounded order. The formulation for each fixed $\eta$ is what permits the exponent two in \eqref{eq:liminf}.
\end{remark}

\begin{proof}[Proof of Corollary~\ref{cor:fixedd}]
For part~(1), let $q$ be a prime power, and take $s=q$, $r=1$ and a fixed $d\ge2$. Then $\ell=(q^d-1)/(q-1)=q^{d-1}+q^{d-2}+\cdots+1\le dq^{d-1}$, so \eqref{eq:lengthcondition} holds with $\varepsilon=1/2$ for all sufficiently large $q$. Moreover, $\ell-1=q^{d-1}(1+O(1/q))$, $(s-1)/(q-1)=1$ and $a_{q,q}=q/\sqrt e-1=(q/\sqrt e)(1+O(1/q))$. Hence,
\[\frac{(\ell-1)(s-1)a_{q,s}}{6(q-1)\ln2}
=\frac{(\ell-1)(q/\sqrt e-1)}{6\ln2}
=\frac{q^d}{6\sqrt e\ln2}\bigl(1+O(1/q)\bigr),\]
and part~(1) follows from Theorem~\ref{thm:general}.
	
For part~(2), fix $r$, let $s\to\infty$ through prime powers, let $q=rs$, and take $d=2$. Then $\ell=s+1\le q+1$, so \eqref{eq:lengthcondition} holds with $\varepsilon=1/2$ for all sufficiently large $q$. Moreover, $\ell-1=s$, $(s-1)/(q-1)=(1+O(1/s))/r$ and $a_{q,s}=se^{-1/(2r)}(1+O(1/s))$. Hence,
\[\frac{(\ell-1)(s-1)a_{q,s}}{6(q-1)\ln2}=\frac{s^2e^{-1/(2r)}}{6r\ln2}\bigl(1+O(1/s)\bigr)
=\left(\frac{e^{-1/(2r)}}{6r^3\ln2}+o(1)\right)q^2,\]
since $s=q/r$. Part~(2) now follows from Theorem~\ref{thm:general}.
\end{proof}

\begin{remark}[The extended Reed--Solomon code]\label{rem:rs}
Let $s=q$ be a prime power and $d=2$. The one-dimensional subspaces of $\F_q^2$ are spanned by $(1,0)$ and by the vectors $(a,1)$ with $a\in\F_q$. With these representatives and $u=(x,y)$, the coordinate forms $u\cdot v$ are $x$ and $y+ax$ for $a\in\F_q$. Since $r=1$, the group $H$ is trivial. Hence, the code \eqref{eq:simplex-code} is
\begin{equation}\label{eq:rscode}
C_{\mathrm{RS}}=\{(x,(y+ax)_{a\in\F_q}):x,y\in\F_q\},\qquad\ell=q+1,
\end{equation}
the extended Reed--Solomon code of dimension two. Since $d=2$, all its triples of coordinates are collinear, and by Corollary~\ref{cor:fixedd}(1) with $d=2$, this inner code yields the increment $(q^2/(6\sqrt e\ln2)+o(q^2))p_q^3$. The same code is the case of a prime power $q$ in Section~\ref{sec:asym}, and Lemma~\ref{lem:bonferroni} recovers this increment directly from the Bonferroni inequality of order four. Section~\ref{sec:counts} evaluates this inner code exactly for small $q$.
\end{remark}

\section{Proof of Theorem~\ref{thm:asym} and Corollary~\ref{cor:special}}\label{sec:asym}

In this section, we prove Theorem~\ref{thm:asym} and Corollary~\ref{cor:special}. The inner code has length $n_q+1\le q+1$, which is short compared with the high-dimensional codes of Section~\ref{sec:simplex}. For such lengths, every triple of coordinates may be dependent, but the Bonferroni inequality of order four suffices. First, we prove a general bound for pairwise uniform codes of length at most $q+1$ in which all triples of coordinates have the same separation probability. Then we construct the product-ring codes and show that, after relabelling the symbols, every projection onto three coordinates becomes a code $C_U(G)$ with $|U|\le2$. Finally, we evaluate these projections by Lemma~\ref{lem:cosettriples} and complete the proof.

\begin{lemma}\label{lem:bonferroni}
Let $\mathcal Q$ be an infinite set of integers, and let $\gamma>0$ be fixed. For each $q\in\mathcal Q$, let $C\subseteq\cA^\ell$ be a pairwise uniform code, where $|\cA|=q$ and the length $\ell$ may depend on $q$ and satisfies $3\le\ell\le q+1$. Suppose that any three distinct coordinates separate $q$ independent uniform words of $C$ with the same probability $t_3$, and that $t_3=(\gamma+o(1))qp_q^3$ as $q\to\infty$ through $\mathcal Q$. Then, as $q\to\infty$ through $\mathcal Q$,
\[R_q\ge B_q+(1+o(1))\,\frac{\gamma(\ell-1)(\ell-2)}{6\ln2}\,p_q^3.\]
\end{lemma}

\begin{proof}
Sample $q$ words of $C$ independently and uniformly, and let $E_j$ be the event that the $j$th coordinate separates them. Write $p=p_q$, and let $S_h$ be as in \eqref{eq:bonferroni-general}. By \eqref{eq:pairwise}, we have $S_1=\ell p$ and $S_2=\binom\ell2p^2$. By hypothesis, $S_3=\binom\ell3t_3$. Let $t_4^*$ be the largest probability of an intersection of four of these events, with $t_4^*=0$ if $\ell=3$. Then $S_4\le\binom\ell4t_4^*$.
	
We first compare $t_4^*$ with $t_3$. By Lemma~\ref{lem:matching} with $h=4$, every intersection of four of these events has probability at most $p\exp(-3q+q\omega_4(q))$, and hence so does $t_4^*$. Moreover, $q\omega_4(q)=O(q^{2/3})$ by \eqref{eq:entropy-error-rate}. For all sufficiently large $q\in\mathcal Q$, we have $t_3\ge\gamma qp^3/2$, and $p^2\ge2\pi qe^{-2q}$ by \eqref{eq:stirling}. Hence,
\[\frac{t_4^*}{t_3}\le\frac{2e^{-3q+O(q^{2/3})}}{\gamma qp^2}\le\exp(-q+o(q)).\]
	
By the Bonferroni inequality \eqref{eq:bonferroni-general} with $K=4$,
\[F(C)\le1-S_1+S_2-S_3+S_4\le1-\ell p+\binom\ell2p^2-\binom\ell3t_3+\binom\ell4t_4^*.\]
By \eqref{eq:binomial-truncation} with $K=3$, the truncated sum $P_3(p)=1-\ell p+\binom\ell2p^2-\binom\ell3p^3$ is at most $(1-p)^\ell$. Writing $\binom\ell3t_3=\binom\ell3p^3+\binom\ell3(t_3-p^3)$, we obtain $F(C)\le(1-p)^\ell-Y$, where $Y=\binom\ell3(t_3-p^3)-\binom\ell4t_4^*.$
	
We now show that $Y=(1+o(1))\binom\ell3\gamma qp^3$. Since $\binom\ell4/\binom\ell3=(\ell-3)/4\le q/4$, it follows from the above bound for $t_4^*/t_3$ that
\[ \frac{\binom\ell4t_4^*}{\binom\ell3t_3} \le\frac q4\exp(-q+o(q))=o(1).\]
Moreover, $p^3/t_3=1/((\gamma+o(1))q)=o(1)$. Hence, $Y=\binom\ell3t_3(1-o(1))=(1+o(1))\binom\ell3\gamma qp^3$. In particular, $Y>0$ for all sufficiently large $q\in\mathcal Q$. Since $\ell\le q+1$, we also have $Y\le\binom\ell3t_3\le(q+1)^3\cdot2\gamma qp^3$, which tends to zero because $q^4p^3\to0$, whereas $(1-p)^\ell\ge1-(q+1)p\to1$. Thus, $0<Y<(1-p)^\ell$ for all sufficiently large $q\in\mathcal Q$. Applying Lemma~\ref{lem:increment}(i) to the bound $F(C)\le(1-p)^\ell-Y$, we obtain
\[ R_q\ge B_q+\frac{Y}{\ell(q-1)\ln2}
=B_q+(1+o(1))\frac{\binom\ell3\gamma q}{\ell(q-1)\ln2}p^3.\]
Since
\[ \frac{\binom\ell3q}{\ell(q-1)}=\frac{(\ell-1)(\ell-2)}6\cdot\frac q{q-1}
=(1+o(1))\frac{(\ell-1)(\ell-2)}6,\]
this proves the lemma.
\end{proof}

We now construct the product-ring inner codes behind Theorem~\ref{thm:asym}. Let $q\ge3$. If $q\not\equiv2\pmod4$, let $W=\{0\}$ and $m=q$. If $q\equiv2\pmod4$, let $W=\Z/2\Z$ and $m=q/2$. In both cases, $m=q/|W|$.
Write $m=q_1\cdots q_k$, where $q_1,\ldots,q_k$ are powers of distinct primes. Then $n_q=\min_iq_i$, by the definition of $n_q$. Moreover, every $q_i$ is odd or divisible by $4$, as noted before Theorem~\ref{thm:asym}. Let $B$ be the commutative ring $\F_{q_1}\times\cdots\times\F_{q_k}$ of order $m$. An element of $B$ is a unit exactly when all its components are nonzero. For each $i$, fix an injection $\varphi_i\colon[n_q]\to\F_{q_i}$, which exists since $n_q\le q_i$, and let $a_j=(\varphi_1(j),\ldots,\varphi_k(j))\in B$ for $j\in[n_q]$. If $j\ne j'$, then every component of $a_j-a_{j'}$ is nonzero, and hence $a_j-a_{j'}$ is a unit of $B$.

On $B^2$, define the linear forms $L_0(x,y)=x$ and $L_j(x,y)=y+a_jx$ for
$j\in[n_q]$. Let $\ell=n_q+1$, and use the alphabet $\cA=W\times B$, which has $q$ elements. Define
\begin{equation}\label{eq:productcode}
C_q=\bigl\{\bigl((\varepsilon_j,L_j(x,y))\bigr)_{0\le j\le n_q}:
x,y\in B,\ \varepsilon_0,\ldots,\varepsilon_{n_q}\in W\bigr\}
\subseteq\cA^{\ell}.
\end{equation}
If $W=\{0\}$, we identify $\cA$ with $B$, and $C_q$ is the $B$-linear code whose coordinates are given by the forms $L_0,\ldots,L_{n_q}$. In particular, if $q$ is a prime power, then $B=\F_q$, the map $\varphi_1$ is a bijection, and $\{a_1,\ldots,a_q\}=\F_q$. In this case, $C_q$ is, up to the order of the coordinates, the extended Reed--Solomon code
\eqref{eq:rscode}.

\begin{lemma}\label{lem:productcode}
Let $G=W\times B$ be the additive group of $\cA$, and let $U=W\times\{0\}$.
\begin{enumerate}
\item [(i)] For any two distinct coordinates, every pair of symbols has exactly $|W|^{\ell-2}$ preimages under the corresponding projection of $C_q$. In particular, $C_q$ is pairwise uniform.
\item [(ii)] For any three distinct coordinates, every word in the image of the corresponding projection has exactly $|W|^{\ell-3}$ preimages. After the symbols in each of the three coordinates are relabelled by a bijection of $\cA$, this image becomes $C_U(G)$.
\end{enumerate}
\end{lemma}

\begin{proof}
We first show that the values of any two of the forms $L_0,\ldots,L_{n_q}$ determine $(x,y)$. For $L_0$ and $L_j$, we have $x=L_0$ and $y=L_j-a_jL_0$. For $L_j$ and $L_{j'}$ with $j\ne j'$, we have $L_j-L_{j'}=(a_j-a_{j'})x$, and since $a_j-a_{j'}$ is a unit, we obtain $x=(a_j-a_{j'})^{-1}(L_j-L_{j'})$ and then $y=L_j-a_jx$. Conversely, all values $(x,y)\in B^2$ occur. Hence, the symbols in two coordinates of a word of $C_q$ determine $x$, $y$ and the two corresponding components in $W$, while the remaining $\ell-2$ components in $W$ are free. This proves (i).
	
Any three of the forms satisfy a linear relation whose coefficients are units of $B$, namely
\[ (a_j-a_{j'})L_0-L_j+L_{j'}=0,\quad(a_{j'}-a_{j''})L_j+(a_{j''}-a_j)L_{j'}+(a_j-a_{j'})L_{j''}=0\]
for distinct $j,j',j''\in[n_q]$. The first relation follows from $L_j-L_{j'}=(a_j-a_{j'})L_0$. In the second relation, the coefficients sum to zero, so the terms involving $y$ cancel, and the coefficient of $x$ is $(a_{j'}-a_{j''})a_j+(a_{j''}-a_j)a_{j'}+(a_j-a_{j'})a_{j''}=0$.
In each of the three coordinates, multiply the $B$-component of every symbol by the corresponding coefficient and keep the $W$-component. This relabels the symbols by bijections of $\cA$, and the relabelled image consists of triples $((e_1,u),(e_2,v),(e_3,w))$ with $u+v+w=0$. Since any two of the three forms determine $(x,y)$, the image has $|W|^3|B|^2$ elements, and so does the set of all such triples. Thus, the relabelled image is the set of all such triples, which is $C_U(G)$. The remaining $\ell-3$ components in $W$ are free, and this proves (ii).
\end{proof}

For $q\not\equiv2\pmod4$, the length $\ell=n_q+1$ is the largest possible for linear codes of this type.

\begin{remark}\label{rem:ringlength}
Let $R$ be a finite commutative ring of order $q$, and let $M$ be a $2\times\ell$ matrix over $R$ such that any two coordinates of $(x,y)M$ determine $(x,y)\in R^2$. Equivalently, all $2\times2$ minors of $M$ are units, since a square matrix over $R$ is invertible exactly when its determinant is a unit. By \cite[Theorem 8.7]{AtiyahMacdonald1969}, $R$ is a direct product of finite local rings, and the order of each factor is a power of the characteristic of its residue field, because the maximal ideal of a finite local ring is nilpotent and the successive quotients of its powers are vector spaces over the residue field. If $\nu(q)=p^a$, then the factors with residue characteristic $p$ have total order $p^a$, since $p^a$ is the exact power of $p$ dividing $q$, so one of them has a residue field of order at most $\nu(q)$. Reducing $M$ modulo the corresponding maximal ideal $\mathfrak m$ of $R$, we obtain a $2\times\ell$ matrix over the field $R/\mathfrak m$ whose $2\times2$ minors are nonzero. Hence, its columns span pairwise distinct one-dimensional subspaces of $(R/\mathfrak m)^2$, and $\ell\le|R/\mathfrak m|+1\le\nu(q)+1$. For $q\not\equiv2\pmod4$, the codes $C_q$ thus have the largest length among such linear codes. For $q\equiv2\pmod4$, we have $\nu(q)=2$, so such codes have length at most $3$. This is why the components in $W$ are needed in that case.
\end{remark}

By Lemma~\ref{lem:productcode}(ii), the three-coordinate intersections of $C_q$ are described by Lemma~\ref{lem:cosettriples} with $G=W\times B$ and $U=W\times\{0\}$. Two other cases of Lemma~\ref{lem:cosettriples} concern the inner codes of Xing and Yuan. For $U=\{0\}$ and $\sigma_G=0$, it is \eqref{eq:zero-triples}. For $q=2m$ with $m$ odd, $G=\Z/q\Z$ and $U=\{0,m\}$, we have $\sigma_G=q(q-1)/2\equiv m\pmod q$, so $\sigma_G\in U$, and $C_U(G)$ is the inner code
\begin{equation}\label{eq:two-cosets}
	C_3=\{(x,y,z)\in(\Z/q\Z)^3:x+y+z\in\{0,m\}\}
\end{equation}
of \cite[Section 5]{XY2023}. Thus, by Lemma~\ref{lem:cosettriples} with $h=2$,
\begin{equation}\label{eq:lift-t3}
	t_3(\Z/q\Z,\{0,m\})=\left(\frac{1}{2e^{1/4}}+o(1)\right)qp_q^3.
\end{equation}
In this case, a function $f\colon[q]\to\{0,m\}$ taking the value $m$ exactly $k$ times has $\sum_if(i)=km$, so it is admissible exactly when $k$ is odd, and then $\mathfrak S_G(f)=\exp(-(k^2+(q-k)^2)/(2q^2))$. The dependence of $\mathfrak S_G(f)$ on $f$ is essential here. 
In \cite[Proposition 4.3]{XY2023}, the count $N_G(f)$ is recorded with the constant $e^{-1/2}$ for every admissible $f$. Since $\sum_{g\in G}|f^{-1}(g)|^2\le q^2$, Theorem~\ref{thm:eberhard} shows that $e^{-1/2}$ is a lower bound for $\mathfrak S_G(f)$, with equality only for constant $f$. Together with the uniformity in \eqref{eq:eberhard-uniform}, this lower bound suffices for \cite[Theorem 5.2]{XY2023}, but it loses the factor $e^{1/4}$ in \eqref{eq:lift-t3}. The code \eqref{eq:two-cosets} is pairwise uniform, since any two of $x,y,z$ are arbitrary and the third lies in a coset of $\{0,m\}$. Already for this code, applying Lemma~\ref{lem:bonferroni} with $\ell=3$ and \eqref{eq:lift-t3}, we obtain
\[ R_q\ge B_q+\left(\frac{1}{6e^{1/4}\ln2}+o(1)\right)p_q^3,\]
which improves \eqref{eq:xy-exceptional} by the factor $e^{1/4}$.

\begin{remark}\label{rem:cosetoptimal}
Let $G$ be an abelian group of order $q$. Every subgroup $U$ with $\sigma_G\in U$ contains the subgroup $\langle\sigma_G\rangle$, which is $\{0\}$ if $\sigma_G=0$ and $\{0,\sigma_G\}$ otherwise. If $U\ne\langle\sigma_G\rangle$, then $|U|\ge2|\langle\sigma_G\rangle|$, and hence
\[ \frac{e^{-1/(2|U|)}}{|U|}<\frac{1}{2|\langle\sigma_G\rangle|}
<\frac{e^{-1/2}}{|\langle\sigma_G\rangle|}
\le\frac{e^{-1/(2|\langle\sigma_G\rangle|)}}{|\langle\sigma_G\rangle|}.\]
Since the error in Lemma~\ref{lem:cosettriples} is uniform over all $G$ and $U$, it follows that, for all sufficiently large $q$, the largest value of $t_3(G,U)$ over subgroups $U\le G$ is attained at $U=\langle\sigma_G\rangle$. Moreover, $\sigma_G$ is the sum of the elements of order at most two in $G$, since the other elements cancel in pairs $\{g,-g\}$. These elements form an elementary abelian $2$-group of some rank $k$, whose elements sum to zero unless $k=1$. Since $k=1$ exactly when the Sylow $2$-subgroup of $G$ is cyclic and nontrivial, we have $\sigma_G\ne0$ exactly in this case. This is always the case if $q\equiv2\pmod4$, while for $q\not\equiv2\pmod4$ the group $G=W\times B$ of Lemma~\ref{lem:productcode} has $\sigma_G=0$. Thus, the constants $\gamma_q$ in Theorem~\ref{thm:asym} are the largest asymptotic values of $t_3(G,U)/(qp_q^3)$ obtainable from Lemma~\ref{lem:cosettriples} over all abelian groups $G$ of order $q$ and all subgroups $U\le G$.
\end{remark}

\begin{proof}[Proof of Theorem~\ref{thm:asym}]
Let $C_q$ be the code \eqref{eq:productcode}, of length $\ell=n_q+1$. Since $3\le n_q\le m\le q$, we have $4\le\ell\le q+1$. By Lemma~\ref{lem:productcode}(i), the code $C_q$ is pairwise uniform. By Lemma~\ref{lem:productcode}(ii), for any three coordinates the projection of a uniform word of $C_q$ is uniform on the image, and after relabelling the symbols, this image is $C_U(G)$ with $G=W\times B$ and $U=W\times\{0\}$. Relabelling the symbols in each coordinate does not change the probability of separation. Hence, any three coordinates of $C_q$ separate $q$ independent uniform words of $C_q$ with probability $t_3(G,U)$.
	
We now compute $\sigma_G$. If $q\not\equiv2\pmod4$, then $W=\{0\}$ and $U=\{0\}$. Each $q_i$ is odd or divisible by $4$, so the elements of each $\F_{q_i}$ sum to zero, and hence $\sigma_G=0$. If $q\equiv2\pmod4$, then $W=\Z/2\Z$ and $|B|=m$ is odd. Each element of $W$ occurs $m$ times as a first component, and each element of $B$ occurs twice as a second component. Hence, $\sigma_G=(m\cdot1,2\sigma_B)=(1,0)\in U$, and $|U|=2$. In both cases, it follows from Lemma~\ref{lem:cosettriples} that $t_3(G,U)=(\gamma_q+o(1))qp_q^3$. Applying Lemma~\ref{lem:bonferroni} separately to the sets of integers $q\not\equiv2\pmod4$ and $q\equiv2\pmod4$, with $\gamma=\gamma_q$ and $(\ell-1)(\ell-2)=n_q(n_q-1)$, we obtain the bound in Theorem~\ref{thm:asym}.
\end{proof}

\begin{proof}[Proof of Corollary~\ref{cor:special}]
Since $n_q\ge3$, we have $n_q(n_q-1)/6\ge1$, and both parts follow from Theorem~\ref{thm:asym}.
\end{proof}

\begin{remark}\label{rem:xy-lemma}
The improvement over \eqref{eq:xy-regular} is already visible for the length-four inner codes of \cite[Theorems 4.7 and 4.8]{XY2023}. These are the codes $C=\{(x,y,x+y,x+cy):x,y\in A\}$, where $A=\Z/q\Z$ and $c=-1$ if $q$ is odd, while $A=\F_{2^a}\times\Z/m'\Z$ and $c=(\alpha,-1)$ with $\alpha\in\F_{2^a}\setminus\F_2$ if $q=2^am'$ with $a\ge2$ and $m'$ odd, the second factor being omitted if $m'=1$. Write $L_1=x$, $L_2=y$, $L_3=x+y$ and $L_4=x+cy$. The four triples of these forms satisfy $L_1+L_2-L_3=0$, $L_1+cL_2-L_4=0$, $(c-1)L_1-cL_3+L_4=0$, $(c-1)L_2+L_3-L_4=0$,
and all coefficients are units because $c$ and $c-1$ are units. As in the proof of Lemma~\ref{lem:productcode}, each triple projection of $C$ is therefore equivalent to $C_{\{0\}}(A)$, and the elements of $A$ sum to zero. Let $a_q$ be the number of $q$-element subsets of $C$ separated by each of three given coordinates, and let $b_q$ be the number separated by all four coordinates. A $q$-element subset separated by one given coordinate is obtained by choosing, for each of the $q$ symbols, one of the $q$ words with this symbol, so there are $q^q$ of them. A subset separated by two given coordinates is the inverse image of the graph of a permutation of $A$, so there are $q!$ of them. By inclusion--exclusion over the four coordinates, we have $|\cS(C)|=4q^q-6q!+4a_q-b_q.$
The estimate in \cite[Lemma 4.6]{XY2023} uses $b_q\le a_q$, and hence $4a_q-b_q\ge3a_q$. By Lemmas~\ref{lem:cosettriples} and~\ref{lem:matching}, we have $b_q/a_q\le\exp(-q+o(q))$, and hence $4a_q-b_q=(4-o(1))a_q$. Thus, applying Lemma~\ref{lem:bonferroni} with $\ell=4$, we obtain the constant $1/(\sqrt e\ln2)$ of Corollary~\ref{cor:special}(1) already for these codes.
For odd $q$, the number $b_q$ counts the configurations of $q$ nonattacking queens on the toroidal $q\times q$ board, and in this case the required upper bound already follows from \cite[Corollary 3.2(1)]{Luria2017}.
\end{remark}

\section{Proof of Theorem~\ref{thm:limit}}\label{sec:limit}

In this section, we prove Theorem~\ref{thm:limit} and deduce its consequences for the length of inner codes with large increments and for MDS codes of dimension two. Since Theorem~\ref{thm:limit} bounds $\Phi(C)$ from above, we need a lower bound for $F(C)$, and we obtain it from the Bonferroni inequality of order three. By pairwise uniformity, the only terms of this inequality that differ from the independent model are the intersections of three separation events, and by Lemma~\ref{lem:matching} with $h=3$, each of them has probability at most $\theta_q$, independently of the pairwise uniform code and of the three coordinates.

\begin{proof}[Proof of Theorem~\ref{thm:limit}]
Let $E_j$ be the event that the $j$th coordinate separates $q$ independent uniform words of $C$, write $p=p_q$, and let $S_h$ be as in \eqref{eq:bonferroni-general}. By \eqref{eq:pairwise}, we have $S_1=\ell p$ and $S_2=\binom\ell2p^2$. By Lemma~\ref{lem:matching} with $h=3$ and by \eqref{eq:omega3}, every intersection of three of these events has probability at most $p\exp(-2q+q\omega_3(q))\le\theta_q$. Hence, $S_3\le\binom\ell3\theta_q$. By the Bonferroni inequality \eqref{eq:bonferroni-general} with $K=3$,
\[ F(C)\ge1-S_1+S_2-S_3\ge1-\ell p+\binom\ell2p^2-\binom\ell3\theta_q.\]
By \eqref{eq:binomial-truncation} with $K=2$, the truncated sum $P_2(p)=1-\ell p+\binom\ell2p^2$ is at least $(1-p)^\ell$. Therefore, $F(C)\ge(1-p)^\ell-Y$, where $Y=\binom\ell3\theta_q$.
By the hypotheses, $Y\le\ell^3\theta_q/6\le1/6$, and by Bernoulli's inequality, $(1-p)^\ell\ge1-\ell p\ge1/2$. Hence, $(1-p)^\ell-Y\ge1/3$, and by Lemma~\ref{lem:increment}(ii),
\[ \Phi(C)\le B_q+\frac{3Y}{\ell(q-1)\ln2} \le B_q+\frac{\ell^2\theta_q}{2(q-1)\ln2}.\]
Since $2\ln2>1$, this proves the theorem.
\end{proof}

\begin{remark}\label{rem:explength}
Let $0<\eta<1$, and suppose that a pairwise uniform inner code $C$ of length $\ell$ yields an increment of at least $\exp(-(2+\eta)q)$ over $B_q$. Then $\ell\ge3$, since $F(C)=(1-p_q)^2$ for $\ell=2$ by inclusion--exclusion and \eqref{eq:pairwise}, and $F(C)\ge1-p_q$ for $\ell=1$. By \eqref{eq:stirling}, we have $\theta_q=\exp(-3q+O(\sqrt q))$. If the hypotheses of Theorem~\ref{thm:limit} fail, then $\ell>1/(2p_q)$ or $\ell>\theta_q^{-1/3}$, and hence $\ell\ge\exp(q-O(\sqrt q))$. Otherwise, it follows from Theorem~\ref{thm:limit} that $\ell^2\theta_q\ge(q-1)\exp(-(2+\eta)q)$, and hence $\ell\ge\exp((1-\eta)q/2-O(\sqrt q))$. Thus, the last bound holds in both cases.
\end{remark}

\begin{remark}\label{rem:limit}
Let $C\subseteq\cA^\ell$ be an MDS code of dimension two with $\ell\ge3$. Since the first two coordinates determine a codeword, each of the other $\ell-2$ coordinates is a function of the first two, and these functions form $\ell-2$ mutually orthogonal Latin squares of order $q$. There are at most $q-1$ such squares \cite{ColbournDinitz2007}, and hence $\ell\le q+1$. Moreover, the projection of $C$ onto any two coordinates is a bijection onto $\cA^2$, so $C$ is pairwise uniform.

The hypotheses of Theorem~\ref{thm:limit} become weaker as $\ell$ decreases, so it suffices to check them for $\ell=q+1$ and $q\ge12$. Since $p_{q+1}/p_q=(1+1/q)^{-q}\le1/2$ by Bernoulli's inequality, the ratio of consecutive values of $(q+1)p_q$ is at most $(q+2)/(2(q+1))<1$. Together with $5p_4=15/32$, this gives $(q+1)p_q<1/2$ for all $q\ge4$. By \eqref{eq:stirling}, we have $\ln((q+1)^3\theta_q)\le g(q)$, where
\[ g(q)=3\ln(q+1)+\tfrac12\ln(2\pi q)+\frac1{12q}-3q+2(3+\ln2)\sqrt q.\]
Since $g'(q)<3/(q+1)+1/(2q)+(3+\ln2)/\sqrt q-3<0$ for $q\ge12$, the function $g$ is decreasing for $q\ge12$, and a direct computation gives $g(12)<-0.54$. Hence, $(q+1)^3\theta_q<1$ for all $q\ge12$. Moreover, $p_q^2\ge2\pi qe^{-2q}$ by \eqref{eq:stirling}, so that $\theta_q\le p_q^3e^{2(3+\ln2)\sqrt q}/(2\pi q)$. Since $\ell\le q+1$ and $(q+1)^2\le2\pi q(q-1)$, it follows from Theorem~\ref{thm:limit} that
\[-\frac{\log_2F(C)}{\ell(q-1)}
\le B_q+\frac{(q+1)^2}{2\pi q(q-1)}e^{2(3+\ln2)\sqrt q}p_q^3
\le B_q+e^{2(3+\ln2)\sqrt q}p_q^3\]
for every MDS code of dimension two and every $q\ge12$.

The projection of such a code onto three coordinates consists of the triples $(x,y,L(x,y))$ for a Latin square $L$ of order $q$, and a $q$-element subset is separated by all three coordinates exactly when it is a transversal of $L$. In this case, Lemma~\ref{lem:matching} with $h=3$ is a quantitative form of the upper bound $((1+o(1))q/e^2)^q$ for the number of transversals of a Latin square of order $q$, due to Taranenko \cite{Taranenko2015}. A simpler proof of this bound is due to Glebov and Luria \cite{GlebovLuria2016}. The factor $e^{O(\sqrt q)}$ above comes only from this bound. For the codes of Theorem~\ref{thm:asym} with $q\not\equiv2\pmod4$, every triple projection has $t_3=(e^{-1/2}+o(1))qp_q^3$, which corresponds to $(e^{-1/2}+o(1))q(q!)^2/q^q$ transversals. An upper bound of the form $q^{O(1)}(q!)^2/q^q$ for the number of transversals of every Latin square of order $q$ would replace $e^{O(\sqrt q)}$ by $q^{O(1)}$ for MDS codes of dimension two.
\end{remark}

\section{Exact counts for extended Reed--Solomon inner codes}\label{sec:counts}

In this section, we compute $|\cS(C_{\mathrm{RS}})|$ exactly for small prime powers $q$, where $C_{\mathrm{RS}}$ is the code \eqref{eq:rscode}. By Lemma~\ref{lem:concatenation} and \eqref{eq:failure-count}, these counts lead to explicit lower bounds for $R_q$, which the asymptotic arguments in Sections~\ref{sec:simplex} and~\ref{sec:asym} do not provide for a particular small alphabet size. Since a direct enumeration of all $\binom{q^2}{q}$ subsets of size $q$ is infeasible except for very small $q$, we count incidences between separated subsets and their separating coordinates in two ways, using the transitive action of $\GL_2(\F_q)$ on coordinates and on pairs of coordinates. This reduces the count to histograms over functions or permutations of $\F_q$ (Lemma~\ref{lem:exactcount}), and a further reduction by affine symmetries makes the enumeration feasible for $q\le13$.

For a function $f\colon\F_q\to\F_q$, set
\[k(f)=1+|\{a\in\F_q:x\mapsto f(x)+ax\text{ is a permutation}\}|.\]
The summand $1$ accounts for the coordinate $x$. Let $H_k$ be the number of functions $f$ with $k(f)=k$, and let $J_k$ be the number of permutations $f$ with $k(f)=k$. Thus, $\sum_kH_k=q^q$ and $\sum_kJ_k=q!$. Moreover, $J_1=0$, since $k(f)\ge2$ for every permutation $f$, the coordinate $y$ corresponding to $a=0$.

\begin{lemma}\label{lem:exactcount}
For the code in \eqref{eq:rscode}, we have
\begin{align}
|\cS(C_{\mathrm{RS}})|&=(q+1)\sum_{k=1}^{q+1}\frac{H_k}{k},\label{eq:count-functions}\\
|\cS(C_{\mathrm{RS}})|&=(q+1)q^q-q(q+1)\sum_{k=2}^{q+1}\frac{J_k}{k}.
\label{eq:count-permutations}
\end{align}
Moreover, the numbers $H_k$ can be computed from $q^{q-2}$ normalized functions, and the numbers $J_k$ from $(q-2)!$ normalized permutations.
\end{lemma}

\begin{proof}
The map $(x,y)\mapsto(x,(y+ax)_{a\in\F_q})$ is a bijection from $\F_q^2$ onto $C_{\mathrm{RS}}$, since the coordinates $x$ and $y$ of a codeword determine it. We identify $q$-element subsets of $C_{\mathrm{RS}}$ with $q$-element subsets of $\F_q^2$ in this way. The coordinates of $C_{\mathrm{RS}}$ correspond to the linear forms $x$ and $y+ax$ with $a\in\F_q$, which represent the $q+1$ one-dimensional subspaces of the dual space of $\F_q^2$. A coordinate separates a subset $P\subseteq\F_q^2$ exactly when the corresponding form is injective on $P$.
	
If the coordinate $x$ separates $P$, then $P$ contains exactly one point with each first coordinate, and hence $P$ is the graph $\{(x,f(x)):x\in\F_q\}$ of a unique function $f$. The coordinate $y+ax$ separates this graph exactly when $x\mapsto f(x)+ax$ is a permutation.
Hence, the number of coordinates separating the graph of $f$ is $k(f)$. Let $S_k$ be the number of $q$-element subsets of $\F_q^2$ separated by exactly $k$ coordinates. 
Let $g\in\GL_2(\F_q)$. For a nonzero linear form $L$, the form $L\circ g$ is nonzero, and hence equals $cL'$ for a unique coordinate form $L'$ and a unique $c\in\F_q^\times$. Since multiplication by $c$ is injective, $L$ separates $g(P)$ exactly when $L\circ g$ separates $P$, that is, exactly when $L'$ separates $P$. Thus, $g$ induces a permutation of the $q+1$ coordinates, and $P$ and $g(P)$ are separated by the same number of coordinates. Since $\GL_2(\F_q)$ acts transitively on ordered bases of the dual space, and hence on ordered pairs of distinct one-dimensional subspaces, the induced action is transitive on coordinates and on pairs of distinct coordinates.
	
We count the pairs $(P,L)$ such that $P$ is separated by exactly $k$ coordinates and $L$ is one of them. Each such $P$ contributes $k$ pairs, so there are $kS_k$ pairs. On the other hand, for each coordinate $L$, the number of subsets $P$ separated by $L$ and by exactly $k$ coordinates in total is the same as for the coordinate $x$, by the transitivity on coordinates. For $L=x$, this number is $H_k$, by the second paragraph. Hence,
\begin{equation}\label{eq:singleincidence}
kS_k=(q+1)H_k.
\end{equation}
Dividing by $k$ and summing over $k\ge1$ proves \eqref{eq:count-functions}.
	
A subset is separated by both coordinates $x$ and $y$ if and only if it is the graph of a permutation $f$, and then it is separated by exactly $k(f)$ coordinates. Counting the pairs $(P,\{L,L'\})$ such that $P$ is separated by exactly $k$ coordinates, including the two distinct coordinates $L$ and $L'$, and using the transitivity on pairs, we obtain
\begin{equation}\label{eq:pairincidence}
\binom{k}{2}S_k=\binom{q+1}{2}J_k\qquad(k\ge2).
\end{equation}
Summing \eqref{eq:singleincidence} over $k$, we obtain $\sum_k kS_k=(q+1)\sum_kH_k=(q+1)q^q$. By \eqref{eq:pairincidence}, $(k-1)S_k=q(q+1)J_k/k$ for $k\ge2$. Hence,
\begin{align*}
|\cS(C_{\mathrm{RS}})|=\sum_{k\ge1}S_k = \sum_{k\ge1}kS_k-\sum_{k\ge2}(k-1)S_k
= (q+1)q^q-q(q+1)\sum_{k\ge2}\frac{J_k}{k}.
\end{align*}
	
For the normalization, let $a,b\in\F_q$. Replacing $f$ by $x\mapsto f(x)+ax+b$ replaces each map $x\mapsto f(x)+a'x$ by $x\mapsto f(x)+(a'+a)x+b$. Since adding a constant does not affect bijectivity, this shifts the set of good directions $a'$ by $-a$ and preserves $k(f)$. The resulting action of the $q^2$ affine functions is free, so every orbit has $q^2$ elements. Every orbit contains a unique function with $f(0)=f(1)=0$, obtained with $b=-f(0)$ and $a=f(0)-f(1)$. Thus, $H_k$ can be computed by enumerating these $q^{q-2}$ functions and multiplying the resulting histogram by $q^2$. For permutations, the transformations $f\mapsto cf+b$ with $c\in\F_q^\times$ and $b\in\F_q$ preserve the permutation property. They also preserve $k(f)$, since $cf+b+a'c\,x=c(f+a'x)+b$ shows that $a'$ is a good direction for $f$ exactly when $ca'$ is one for $cf+b$. This action is free, so each orbit has $q(q-1)$ elements, and each orbit contains a unique permutation with $f(0)=0$ and $f(1)=1$. Enumerating these $(q-2)!$ permutations and multiplying the histogram by $q(q-1)$, we obtain $J_k$.
\end{proof}

For a given function or permutation, testing all directions $a\in\F_q$ takes $O(q^2)$ field operations. Thus, the two enumeration procedures use $O(q^q)$ and $O(q^2(q-2)!)=O(q!)$ field operations, respectively, and $O(q^2)$ working storage. The bit cost of the integer arithmetic and of the arbitrary-precision accumulation is not included in these counts.
For example, for $q=8$ the nonzero values of $H_k$ are $H_1=16461312$, $H_2=312704$, $H_4=3136$, and $H_8=64.$
For $q=9$, they are $H_1=384223824$, $H_2=3135672$, $H_3=54432$, $H_4=5832$, $H_6=648$, and $H_9=81.$
For the numerical evaluation, we compute the numbers $J_k$ from the $(q-2)!$ normalized permutations of Lemma~\ref{lem:exactcount} and substitute them into \eqref{eq:count-permutations}.
By \eqref{eq:failure-count} and Lemma~\ref{lem:concatenation}, we obtain $R_q\ge\underline R_q^{\mathrm{RS}}$, where
\[\underline R_q^{\mathrm{RS}}=-\frac{1}{(q+1)(q-1)}
\log_2\left(1-\frac{q!\,|\cS(C_{\mathrm{RS}})|}{q^{2q}}\right).\]
In Table~\ref{tab:counts}, we list the exact numbers of separated subsets of the code \eqref{eq:rscode} and the corresponding lower bounds $\underline R_q^{\mathrm{RS}}$ on $R_q$. For comparison, it also lists the probabilistic bound $B_q$ and the lower bounds of Xing and Yuan \cite{XY2023}. 

\begin{table}[t]
\centering
\setlength{\belowcaptionskip}{3pt}
\caption{Exact counts and rate bounds.}
\label{tab:counts}
\small
\renewcommand{\arraystretch}{1.25}
\begin{tabular*}{\textwidth}{@{\extracolsep{\fill}}rrlll@{}}
\toprule
\multicolumn{1}{@{}c}{$q$}
& \multicolumn{1}{c}{$|\cS(C_{\mathrm{RS}})|$}
& \multicolumn{1}{c}{$B_q$}
& \multicolumn{1}{c}{Xing--Yuan \cite{XY2023}}
& \multicolumn{1}{c@{}}{$\underline R_q^{\mathrm{RS}}$}\\
\midrule
3  & 84        & 0.1812850396923541 & $0.2119992266387375^{a}$ & 0.2119992266387375\\
4  & 1100      & 0.0473396682908092 & $0.0495860874315075^{a}$ & 0.0495860874315075\\
5  & 17130     & 0.0141227997095659 & $0.0145266449305914^{b}$ & 0.0142073683923850\\
7  & 6453552   & 0.0014760459085397 & $0.0014833905220379^{b}$ & 0.0014771283199558\\
8  & 149566104   & 0.0004959061690047 & $0.000495909^{c}$ & 0.0004959550085065\\
9  & 3858113790       & 0.0001689929050610 & $0.0001689931^{c}$ & 0.0001690016093551\\
11 & 3421113322836    & 0.0000201855739361 & $0.0000201855746^{c}$ & 0.0000201856126515\\
13 & 4239685205859376 & 0.0000024718066359 & \multicolumn{1}{c}{---} & 0.0000024718068259\\
\bottomrule
\end{tabular*}
	
\smallskip
\parbox{\textwidth}{\footnotesize Values with sixteen decimal places are rounded down. $^{a}$Obtained from the same count of separated subsets. $^{b}$Computed from the counts in \cite[Lemma 6.1]{XY2023} for the nonlinear construction of \cite[Section 6]{XY2023} and Lemma~\ref{lem:concatenation}. The values stated in \cite[Corollary 6.2]{XY2023} are $0.01452$ and $0.001483$. $^{c}$As reported in \cite[Table 3]{XY2023}.}
\end{table}

For $q=3,4$, our counts agree with those of \cite{XY2023}, and for $q=3$, the resulting bound is exactly the bound \eqref{eq:KMbound} of K\"orner and Marton. For $q=5,7$, the nonlinear construction in \cite[Section 6]{XY2023} is better than the extended Reed--Solomon inner code. Thus, increasing the length of the MDS inner code does not improve the best available construction for every small alphabet size. For $q=8,9,11$, the bounds $\underline R_q^{\mathrm{RS}}$ improve those reported in \cite[Table 3]{XY2023}. The values $\underline R_q^{\mathrm{RS}}$ are lower bounds obtained from one particular inner code, and they do not determine $R_q$.

\section{Concluding remarks}\label{sec:remaining}

In Theorem~\ref{thm:allq}, we proved that the classical probabilistic lower bound for perfect $q$-hash codes is sharp if and only if $q=2$, which settles the first assertion of the Xing--Yuan conjecture with an explicit bound for every $q\ge3$. In Theorem~\ref{thm:exponential}, using high-dimensional simplex inner codes, we showed that $R_q-B_q\ge\exp(-(2+\eta)q)$ for each fixed $\eta>0$ and all sufficiently large $q$, whereas the previously known increments are of order $\exp(-3q+O(\ln q))$. In Theorem~\ref{thm:asym}, we used MacNeish's product construction to strengthen the bounds obtained from inner codes of length at most $q+1$. In Theorem~\ref{thm:limit}, we bounded the increments obtainable from pairwise uniform inner codes and deduced that such codes need exponential length to reach the scale of Theorem~\ref{thm:exponential}. For polynomial lengths, these increments are at most $e^{O(\sqrt q)}p_q^3$, so the constructions of Theorem~\ref{thm:asym} are optimal for such lengths up to this factor.

The second assertion of the Xing--Yuan conjecture, which concerns the short inner codes specified in \cite[Theorems 4.7, 4.8 and 5.2]{XY2023}, remains open. For the length-four MDS code, with $a_q$ and $b_q$ as in Remark~\ref{rem:xy-lemma}, it amounts to the inequality $4a_q-b_q>\frac{4(q!)^2}{q^q}-\frac{(q!)^3}{q^{2q}}$,
which Lemmas~\ref{lem:cosettriples} and~\ref{lem:matching} establish only for all sufficiently large $q\not\equiv2\pmod4$. For the two-coset inner code \eqref{eq:two-cosets}, the corresponding condition $t_3>p_q^3$ is likewise known only asymptotically, by \eqref{eq:lift-t3}.

It is natural to ask how large the increment $R_q-B_q$ can be. Since $R_q\le q!/q^{q-1}=qp_q$ by the Fredman--Koml\'os bound \cite{FK1984}, we have $R_q-B_q\le e^{-q+O(\ln q)}$ by \eqref{eq:stirling}. Together with \eqref{eq:liminf}, this shows that all limit points of $q^{-1}\ln(R_q-B_q)$ lie in the interval $[-2,-1]$.

\begin{problem}\label{prob:exponent}
Determine the limit points of $q^{-1}\ln(R_q-B_q)$ as $q\to\infty$. In particular, is $\liminf_{q\to\infty}q^{-1}\ln(R_q-B_q)>-2$?
\end{problem}

Within the concatenation method, Remark~\ref{rem:explength} shows that, for $0<\eta<1$, a pairwise uniform inner code with an increment of at least $\exp(-(2+\eta)q)$ has length at least $\exp((1-\eta)q/2-O(\sqrt q))$, whereas our simplex codes have length close to $\exp((1-\eta/2)q)$. Inner lengths beyond $\exp(q)$ would require control of intersections of unbounded order, which our fixed-order Bonferroni argument does not provide.

For MDS codes of dimension two, the factor $e^{O(\sqrt q)}$ in Remark~\ref{rem:limit} comes only from the upper bound for the number of transversals of Latin squares. By the results of Taranenko \cite{Taranenko2015} and of Glebov and Luria \cite{GlebovLuria2016}, the maximum number of transversals of a Latin square of order $q$ is $((1+o(1))q/e^2)^q=e^{o(q)}(q!)^2/q^q$. For comparison, if $G$ is an abelian group of order $q$ with $\sigma_G=0$, then the Cayley table of $G$ has $N_G(0)/q!$ transversals, which is $(e^{-1/2}+o(1))q(q!)^2/q^q$ by \eqref{eq:zero-triples}.

\begin{problem}\label{prob:transversals}
Is it true that every Latin square of order $q$ has at most $q^{O(1)}(q!)^2/q^q$ transversals?
\end{problem}

An affirmative answer would replace the factor $e^{O(\sqrt q)}$ in Remark~\ref{rem:limit} by $q^{O(1)}$.

Finally, two questions concern effectivity and explicitness. Theorems~\ref{thm:spectral} and~\ref{thm:length} are explicit, but the thresholds in Theorems~\ref{thm:exponential}, \ref{thm:general} and~\ref{thm:asym} depend on the asymptotic count of additive triples of bijections in Theorem~\ref{thm:eberhard}, and effective versions would require a quantitative form of that theorem. Moreover, our outer codes are chosen at random, and constructing explicit outer codes that achieve the same rates with these inner codes remains open.

\paragraph{Acknowledgments.} The work of Shitao Li was supported by the National Natural Science Foundation of China under Grants 12526612 and 126011020. The work of Minjia Shi was supported by the National Natural Science Foundation of China under Grants 12471490 and W2611002.
During the preparation and revision of this manuscript, Claude (Anthropic) suggested the directions that led to Theorems~\ref{thm:asym} and~\ref{thm:limit}, namely to extend the short inner codes of Xing and Yuan by MacNeish's product construction, and to bound the increments obtainable from short inner codes such as MDS codes of dimension two. It also assisted in developing the corresponding proofs. The authors take full responsibility for all content of the manuscript, including the mathematical statements, proofs, originality claims, and references.


\begin{thebibliography}{99}
	
\bibitem{ABDN2024}
N. Alon, A. Bishnoi, S. Das and A. Neri,
Strong blocking sets and minimal codes from expander graphs,
\emph{Trans. Amer. Math. Soc.} \textbf{377} (2024), no. 8, 5389--5410.
	
\bibitem{AYZ1995}
N. Alon, R. Yuster and U. Zwick,
Color-coding,
\emph{J. ACM} \textbf{42} (1995), no. 4, 844--856.
	
\bibitem{AtiyahMacdonald1969}
M. F. Atiyah and I. G. Macdonald,
\emph{Introduction to Commutative Algebra},
Addison-Wesley, Reading, MA, 1969.
	
\bibitem{BK2025}
S. Bhandari and A. Khetan,
Improved upper bound for the size of a trifferent code,
\emph{Combinatorica} \textbf{45} (2025), Paper No. 2.
	
\bibitem{BDGP2024}
A. Bishnoi, J. D'haeseleer, D. Gijswijt and A. Potukuchi,
Blocking sets, minimal codes and trifferent codes,
\emph{J. Lond. Math. Soc.} \textbf{109} (2024), no. 6, e12938.
	
\bibitem{BT2026}
A. Bishnoi and I. Tomon,
Explicit constructions of optimal blocking sets and minimal codes,
\emph{Combinatorica} \textbf{46} (2026), Paper No. 13.
	
\bibitem{ColbournDinitz2007}
C. J. Colbourn and J. H. Dinitz (eds.),
\emph{Handbook of Combinatorial Designs}, 2nd ed.,
Chapman \& Hall/CRC, Boca Raton, FL, 2007.
	
\bibitem{CD2021}
S. Costa and M. Dalai,
New bounds for perfect $k$-hashing,
\emph{Discrete Appl. Math.} \textbf{289} (2021),
374--382.
	
\bibitem{DGR2020}
M. Dalai, V. Guruswami and J. Radhakrishnan,
An improved bound on the zero-error list-decoding capacity of the $4/3$ channel,
\emph{IEEE Trans. Inform. Theory} \textbf{66} (2020), no. 2, 749--756.
	
\bibitem{DCD2022}
S. Della Fiore, S. Costa and M. Dalai,
Improved bounds for $(b,k)$-hashing,
\emph{IEEE Trans. Inform. Theory} \textbf{68} (2022), no. 8, 4983--4997.
	
\bibitem{Eberhard2017}
S. Eberhard,
More on additive triples of bijections, arXiv:1704.02407 (2017).
	
\bibitem{EMM2019}
S. Eberhard, F. Manners and R. Mrazovi\'c,
Additive triples of bijections, or the toroidal semiqueens problem,
\emph{J. Eur. Math. Soc.} \textbf{21} (2019), no. 2, 441--463.
	
\bibitem{EMM2022}
S. Eberhard, F. Manners and R. Mrazovi\'c,
An asymptotic for the Hall--Paige conjecture,
\emph{Adv. Math.} \textbf{404} (2022), 108423.
	
\bibitem{FK1984}
M. L. Fredman and J. Koml\'os,
On the size of separating systems and families of perfect hash functions,
\emph{SIAM J. Algebraic Discrete Methods}
\textbf{5} (1984), no. 1, 61--68.
	
\bibitem{GalambosSimonelli1996}
J. Galambos and I. Simonelli,
\emph{Bonferroni-type Inequalities with Applications},
Springer-Verlag, New York, 1996.
	
\bibitem{GlebovLuria2016}
R. Glebov and Z. Luria,
On the maximum number of Latin transversals,
\emph{J. Combin. Theory Ser. A} \textbf{141} (2016), 136--146.
	
\bibitem{GR2022}
V. Guruswami and A. Riazanov,
Beating Fredman--Koml\'os for perfect $k$-hashing,
\emph{J. Combin. Theory Ser. A}
\textbf{188} (2022), 105580.
	
\bibitem{HSG2026}
X. Han, Y. Sun and G. Ge,
Polynomially improved lower bounds for trifferent codes via locally sparse $3$-uniform hypergraphs,
arXiv:2607.26376v2 (2026).
	
\bibitem{Korner1986}
J. K\"orner,
Fredman--Koml\'os bounds and information theory,
\emph{SIAM J. Algebraic Discrete Methods}
\textbf{7} (1986), no. 4, 560--570.
	
\bibitem{KM1988}
J. K\"orner and K. Marton,
New bounds for perfect hashing via information theory,
\emph{European J. Combin.} \textbf{9} (1988), no. 6, 523--530.
	
\bibitem{LPW2009}
D. A. Levin, Y. Peres and E. L. Wilmer,
\emph{Markov Chains and Mixing Times},
American Mathematical Society, Providence, RI, 2009.
	
\bibitem{Luria2017}
Z. Luria,
New bounds on the number of $n$-queens configurations,
arXiv:1705.05225 (2017).
	
\bibitem{MacNeish1922}
H. F. MacNeish,
Euler squares,
\emph{Ann. of Math. (2)} \textbf{23} (1922), no. 3, 221--227.
	
\bibitem{MP2025}
A. M\"uyesser and A. Pokrovskiy,
A random Hall--Paige conjecture,
\emph{Invent. Math.} \textbf{240} (2025), 779--867.
	
\bibitem{Nilli1994}
A. Nilli,
Perfect hashing and probability,
\emph{Combin. Probab. Comput.} \textbf{3} (1994), no. 3, 407--409.
	
	
\bibitem{Robbins1955}
H. Robbins,
A remark on Stirling's formula,
\emph{Amer. Math. Monthly} \textbf{62} (1955), no. 1, 26--29.
	
\bibitem{TaoVu2006}
T. Tao and V. H. Vu,
\emph{Additive Combinatorics},
Cambridge Stud. Adv. Math., vol. 105,
Cambridge University Press, Cambridge, 2006.
	
	
\bibitem{Taranenko2015}
A. A. Taranenko,
Multidimensional permanents and an upper bound on the number of transversals in Latin squares,
\emph{J. Combin. Des.} \textbf{23} (2015), no. 7, 305--320.
	
\bibitem{Taranenko2018}
A. A. Taranenko,
Transversals, plexes, and multiplexes in iterated quasigroups,
\emph{Electron. J. Combin.} \textbf{25} (2018), no. 4, Paper No. 4.30, 17 pp.
	
\bibitem{Taranenko2021}
A. A. Taranenko,
Transversals, near transversals, and diagonals in iterated groups and quasigroups,
\emph{Electron. J. Combin.} \textbf{28} (2021), no. 3, Paper No. 3.48, 22 pp.
	
\bibitem{WX2001}
H. Wang and C. Xing, Explicit constructions of perfect hash families from algebraic curves over finite fields,
\emph{J. Combin. Theory Ser. A} \textbf{93} (2001), no. 1, 112--124.
	
\bibitem{XY2023}
C. Xing and C. Yuan,
Beating the probabilistic lower bound on $q$-perfect hashing,
\emph{Combinatorica} \textbf{43} (2023), 347--366.
	
\end{thebibliography}
\end{document}